\documentclass[letterpaper,12pt]{article}

\usepackage{algorithm}
\usepackage{algpseudocode}
\usepackage[margin=1in]{geometry}  
\usepackage{xspace,exscale,relsize}
\usepackage{fancybox,shadow}
\usepackage{graphicx}
\usepackage{color}
\usepackage{amsthm,amsfonts}
\usepackage{amssymb}
\usepackage{authblk}
\usepackage[title]{appendix}

\usepackage{extarrows}

\usepackage{amsmath,mathrsfs}
\makeatletter
\let\over=\@@over \let\overwithdelims=\@@overwithdelims
\let\atop=\@@atop \let\atopwithdelims=\@@atopwithdelims
\let\above=\@@above \let\abovewithdelims=\@@abovewithdelims
\makeatother
\usepackage{rotating}

	\usepackage{ifpdf}

	\usepackage{subfigure}
	\usepackage{psfrag}
	
	\usepackage{prettyref}
	\usepackage{enumitem}

	\usepackage{tikz}
	\usetikzlibrary{arrows}
	\tikzstyle{int}=[draw, fill=blue!20, minimum size=2em]
	\tikzstyle{dot}=[circle, draw, fill=blue!20, minimum size=2em]
	\tikzstyle{init} = [pin edge={to-,thin,black}]

	\usepackage[
	CJKbookmarks=true,
	bookmarksnumbered=true,
	bookmarksopen=true,
	colorlinks=true,
	citecolor=blue,
	linkcolor=blue,
	anchorcolor=red,
	urlcolor=blue,
	pdfauthor={Soham Jana}
	]{hyperref}
	
	\usepackage[all]{xy}
	
	\usepackage{mathtools}
	
	\usepackage{ifthen}
	\newboolean{aos}
	\setboolean{aos}{FALSE}
	
	\newcommand{\vY}{\boldsymbol{Y}}
	\newcommand{\vx}{\boldsymbol{x}}
	\newcommand{\vy}{\boldsymbol{y}}
	
	\newcommand{\ve}{\boldsymbol{e}}
	\newcommand{\vt}{\boldsymbol{t}}
	
	\newcommand{\dist}{{\sf dist}}
	\newcommand{\vtheta}{\boldsymbol{\theta}}

	\ifx\eqref\undefined
	\newcommand{\eqref}[1]{~(\ref{#1})}
	\fi
	\ifx\mod\undefined
	\def\mod{\mathop{\rm mod}}
	\fi
	
	\usepackage{bm}
	\def\vect#1{\bm{#1}}

	\def\argmin{\mathop{\rm argmin}}
	\def\argmax{\mathop{\rm argmax}}

	\def\EE{\Expect}

	\def\Var{\mathrm{Var}}

	\def\PP{\mathbb{P}}

	\def\eqdef{\triangleq}

	\def\simiid{\stackrel{iid}{\sim}}
	
	\newcommand{\hatGn}{\hat G}

	\newcommand{\abs}[1]{\left| #1 \right|}

	\newcommand{\stepa}[1]{\overset{\rm (a)}{#1}}
	\newcommand{\stepb}[1]{\overset{\rm (b)}{#1}}

	\newcommand{\btheta}{\vect{\theta}}
	\newcommand{\Poi}{\mathrm{Poi}}

	\newcommand{\Unif}{\mathrm{Uniform}}

	\newcommand{\ceil}[1]{{\left\lceil {#1} \right \rceil}}

	\newcommand{\reals}{\mathbb{R}}
	
	\newcommand{\integers}{\mathbb{Z}}
	
	\newcommand{\Expect}{\mathbb{E}}

	\newcommand{\TV}{{\rm TV}}

	\newcommand{\iid}{i.i.d.\xspace}
	\newcommand{\ind}{ind.\xspace}

	\newcommand{\pth}[1]{\left( #1 \right)}
	\newcommand{\qth}[1]{\left[ #1 \right]}
	\newcommand{\sth}[1]{\left\{ #1 \right\}}

	\newcommand{\iiddistr}{{\stackrel{\text{\iid}}{\sim}}}
	\newcommand{\inddistr}{{\stackrel{\text{\ind}}{\sim}}}

	\newcommand{\indc}[1]{{\mathbf{1}_{\left\{{#1}\right\}}}}

	\definecolor{myblue}{rgb}{.8, .8, 1}
	\definecolor{mathblue}{rgb}{0.2472, 0.24, 0.6} 
	\definecolor{mathred}{rgb}{0.6, 0.24, 0.442893}
	\definecolor{mathyellow}{rgb}{0.6, 0.547014, 0.24}

	\newcommand{\calG}{{\mathcal{G}}}

	\newcommand{\calP}{{\mathcal{P}}}

	\newcommand{\mmse}{\text{mmse}}
	\newcommand{\subexpo}{\mathsf{SubE}}
	
	\newrefformat{eq}{(\ref{#1})}
	\newrefformat{thm}{Theorem~\ref{#1}}
	\newrefformat{th}{Theorem~\ref{#1}}
	\newrefformat{chap}{Chapter~\ref{#1}}
	\newrefformat{sec}{Section~\ref{#1}}
	\newrefformat{seca}{Section~\ref{#1}}
	\newrefformat{algo}{Algorithm~\ref{#1}}
	\newrefformat{fig}{Fig.~\ref{#1}}
	\newrefformat{tab}{Table~\ref{#1}}
	\newrefformat{rmk}{Remark~\ref{#1}}
	\newrefformat{clm}{Claim~\ref{#1}}
	\newrefformat{def}{Definition~\ref{#1}}
	\newrefformat{cor}{Corollary~\ref{#1}}
	\newrefformat{lmm}{Lemma~\ref{#1}}
	\newrefformat{propo}{Proposition~\ref{#1}}
	\newrefformat{proper}{Property~\ref{#1}}
	\newrefformat{app}{Appendix~\ref{#1}}
	\newrefformat{apx}{Appendix~\ref{#1}}
	\newrefformat{ex}{Example~\ref{#1}}
	\newrefformat{exer}{Exercise~\ref{#1}}
	\newrefformat{soln}{Solution~\ref{#1}}
	\newrefformat{pt}{Assumption~\ref{#1}}
	
	\def\unifto{\mathop{{\mskip 3mu plus 2mu minus 1mu%
				\setbox0=\hbox{$\mathchar"3221$}%
				\raise.6ex\copy0\kern-\wd0%
				\lower0.5ex\hbox{$\mathchar"3221$}}\mskip 3mu plus 2mu minus 1mu}}

	\ifx\lesssim\undefined
	\def\simleq{{{\mskip 3mu plus 2mu minus 1mu%
				\setbox0=\hbox{$\mathchar"013C$}%
				\raise.2ex\copy0\kern-\wd0%
				\lower0.9ex\hbox{$\mathchar"0218$}}\mskip 3mu plus 2mu minus 1mu}}
	\else
	\def\simleq{\lesssim}
	\fi
	
	\ifx\gtrsim\undefined
	\def\simgeq{{{\mskip 3mu plus 2mu minus 1mu%
				\setbox0=\hbox{$\mathchar"013E$}%
				\raise.2ex\copy0\kern-\wd0%
				\lower0.9ex\hbox{$\mathchar"0218$}}\mskip 3mu plus 2mu minus 1mu}}
	\else
	\def\simgeq{\gtrsim}
	\fi

		\newtheorem{theorem}{Theorem}
		\newtheorem{lemma}[theorem]{Lemma}

		\theoremstyle{definition}

		\newtheorem{remark}{Remark}
		\newtheorem{assumption}{Assumption}
		
		\newif\ifmapx
		{\catcode`/=0 \catcode`\\=12/gdef/mkillslash\#1{#1}}
		\edef\jobnametmp{\expandafter\string\csname ic_apx\endcsname}
		\edef\jobnameapx{\expandafter\mkillslash\jobnametmp}
		\edef\jobnameexpand{\jobname}
		\ifx\jobnameexpand\jobnameapx
		\mapxtrue
		\else
		\mapxfalse
		\fi

		\newcommand{\Regret}{\mathsf{Regret}}

		\renewcommand{\hat}{\widehat}
		\renewcommand{\tilde}{\widetilde}
		
		\newcommand\sfem{\mathsf{emp}}

\begin{document}
	\ifpdf
	\DeclareGraphicsExtensions{.pgf,.jpg}
	\graphicspath{{figures/}{}}
	\fi
	\title{Optimal empirical Bayes estimation for the Poisson model via minimum-distance methods}
	
	\author{Soham Jana, Yury Polyanskiy and Yihong Wu\thanks{
			S.J. is with the Department of ACMS, University of Notre Dame, Notre Dame, IN, email: \url{sjana2@nd.edu}.
			Y.P. is with the Department of EECS, MIT, Cambridge,
			MA, email: \url{yp@mit.edu}. Y.W. is with the Department of Statistics and Data Science, Yale
			University, New Haven, CT, email: \url{yihong.wu@yale.edu}. Y.~Polyanskiy is
			supported in part by the MIT-IBM Watson AI Lab, and the NSF Grants CCF-1717842, CCF-2131115. Y.~Wu is supported in part by the NSF
			Grant CCF-1900507, NSF CAREER award CCF-1651588, and an Alfred Sloan fellowship. }}
	
	\maketitle
	
	\begin{abstract}
		The Robbins estimator is the most iconic and widely used procedure in the empirical Bayes literature for the Poisson model. 
		On one hand, this method has been recently shown to be minimax optimal in terms of the regret (excess risk over the Bayesian oracle that knows the true prior) for various nonparametric classes of priors. On the other hand, it has been long recognized in practice that the Robbins estimator lacks the desired smoothness 
		and monotonicity of Bayes estimators and can be easily derailed by those data points that were rarely observed before.
		Based on the minimum-distance distance method, we propose a suite of empirical Bayes estimators, including 
		the classical nonparametric maximum likelihood, that outperform the Robbins method in a variety of synthetic and real data sets and retain its optimality in terms of minimax regret.
	\end{abstract}

	\maketitle

	\noindent \textbf{Keywords:} {Mixture modeling; Robbins method; Poisson mean estimation; Nonparametric estimation; NPMLE.}

	\section{Introduction}
	\label{sec:EB_introduction}
	Consider the Poisson mean estimation problem. Given observations $Y^n\eqdef (Y_1,\dots,Y_n)$, independently distributed according to the Poisson distribution with mean parameters $\theta^n\eqdef (\theta_1,\dots,\theta_n)$, the goal is to {  estimate} the parameter vector under the squared error loss. It is well known in the literature that the minimax estimator need not be the best choice in practice, unless the observations are known to be generated according to the least favorable prior distribution on the parameter space. A class of shrinkage-type alternative estimators was proposed in the seminal paper of \cite{Rob51,Rob56}, namely the empirical Bayes (EB) methodology. In the regular Bayes setup, which also produces estimators with shrinkage properties, one assumes that the parameter values are independently distributed according to a prior distribution $G$. Then the best estimator under the squared error loss (i.e., the Bayes estimator) of $\theta_j$ is given by the posterior mean $\hat\theta_G(Y_j)=\EE_G\qth{\theta_j|Y_j}$. The EB theory proposes to bypass the assumed knowledge about $G$, which might be unavailable in practice, by approximating the $G$ dependent expressions using the observations.
	For example, in the Poisson model, given a prior distribution $G$ on $\theta$, the posterior mean is of the form 
	\begin{align}\label{eq:EB_bayes_est}
		\hat\theta_{G}(y)
		=\EE_{Y\sim \Poi(\theta),\theta\sim G}\qth{\theta|Y=y}={\int \theta e^{-\theta}{\theta^{y}\over y!}G(d\theta)\over \int e^{-\theta}{\theta^{y}\over y!}G(d\theta)}
		=(y+1){f_G(y+1)\over f_G(y)}.
	\end{align}
	Here $\Poi(\theta)$ denotes the Poisson distribution with mean $\theta$ and marginal density of $Y_j$ is given by
	\begin{align}\label{eq:EB_mixture}
		f_G(y)=\int f_\theta(y)G(d\theta),\quad f_\theta(y)=e^{-\theta}{\theta^y\over y!},y\in\integers_+\eqdef\sth{0,1,\dots}
	\end{align}
	Then, in the EB methodology, one can approximate either $G$ or $f_G$ from the data and plug it into the above formula. The significant achievement of the EB theory is that when the number
	of independent observations is large, it is possible to “borrow strength” from these independent
	(and seemingly unrelated) observations to achieve the asymptotically optimal Bayes risk per coordinate. Since its conception, the theory and methodology of 
	EB has been well developed and widely applied in large-scale data analysis in practice, cf.~e.g.~\cite{efron2001empirical,ver1996parametric,brown2008season,persaud2010comparison}. We refer 
	the reader to the surveys and monographs on the theory and practice of EB \cite{morris1983parametric,casella1985introduction,zhang2003compound,E14,maritz2018empirical,E21}.

\begin{figure}[t]
		\centering
		\includegraphics[height=9cm]{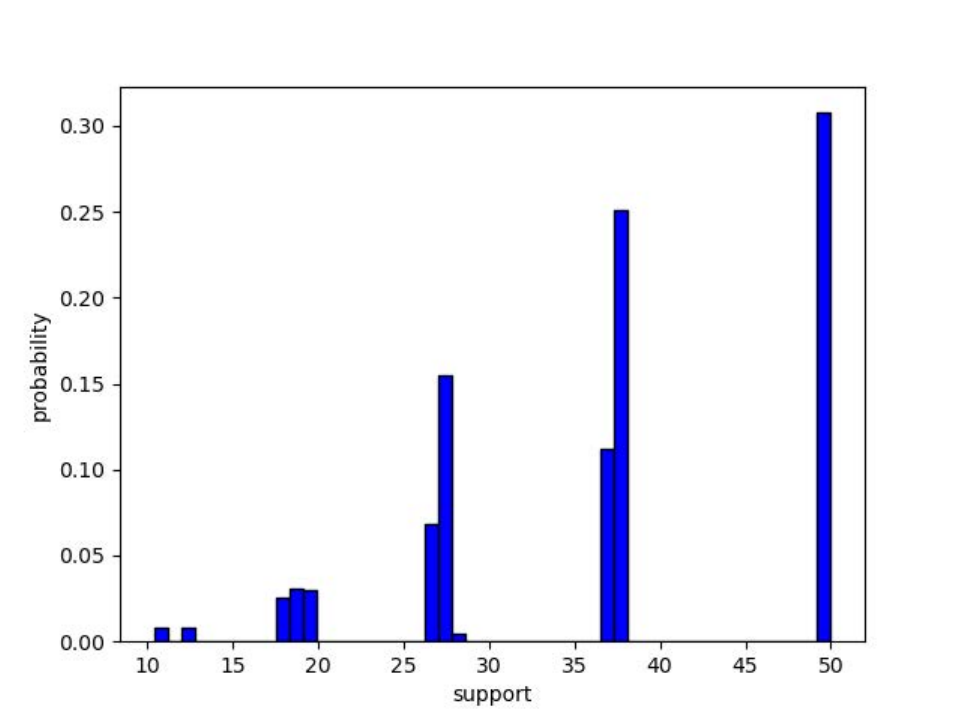}
		\caption{Least favorable prior on [0,50]}
		\label{fig:least_fav_prior}
	\end{figure}

	In particular, to motivate the use of EB methodology in the Poisson settings, we present a real data example where we produce three EB estimators that beat the minimax optimal estimator. We analyze the data on the total number of goals scored in the National Hockey League for the seasons 2017-18 and 2018-19 (the data is available at \href{https://www.hockey-reference.com/}{https://www.hockey-reference.com/}). We consider the statistics of $n=745$ players for whom the data were collected for both seasons. Let $Y_i$ be the total number of goals scored by the $i^{\text{th}}$ player in the season 2017-18. We model $Y_i$ as independently distributed $\Poi(\theta_i)$ random variables, where $\theta_i$'s are independently distributed according to some prior $G$ on $\reals_+$. Based on the observations, we intend to predict the goals scored by each player in the 2018-19 season. Let us explain how \textit{estimation} of $\theta$'s can be used to make \textit{future predictions}. If we assume that player $i$ scores $Z_i \sim \Poi(\theta_i)$ goals in the future year, then predictor $\hat Z_i$ that minimizes mean-square error  (MSE) is clearly $\EE[Z_i | Y_i] = \EE[\theta_i | Y_i] = \hat \theta_G(Y_i)$. Thus, since the EB estimator attempts to approximate $\hat \theta_G$ it can also be used as an estimator for $Z_i$. Note that if the prediction metric is mean \textit{absolute} error (MAE), then the optimal predictor would be a posterior median (under $G$). We do not discuss in this paper EB methods for estimating posterior median and simply reuse the estimator of $\theta_i$ for MAE as well.
	As the number of goals $Y_i$-s in the data are all below 50, for the sake of computation, we can assume that the parameters $\theta_i$-s are supported on $[0,50]$. The minimax estimator in the above Poisson settings with a squared error loss is given by the posterior mean for the least favorable prior (see \prettyref{app:least-favorable} for proof).
	However, even though the minimax estimator is designed to perform optimally for the worst-case scenario, its average performance, particularly with real data, can be overly conservative. Hence, the main advantage of EB is that, in instances far from being the least favorable, one can typically outperform the minimax estimator by being Bayesian with a prior learned from data. Therefore, we will present a comparative study of our EB methods against the minimax estimator and show that in most of our numerical examples, including the hockey-data experiments, there are significant gains in using the EB methodology. To find the least favorable prior, we solve 
		\begin{align}\label{eq:least-favorable-prior-calc}
			\argmax_{G\in \calP([0,h])} \EE[|\theta - \EE[\theta|Y]|^2] = \argmax_{G\in \calP([0,h])} \qth{\EE_G\qth{\theta^2} - \sum_{y=0}^\infty (y+1)^2{(f_G(y+1))^2\over f_G(y)}}
		\end{align} 
		where the maximization is over $G\in \calP([0,h])$, all priors supported on $[0,h]$. For computation, we pick $h=50$, divide parameter space $[0,50]$ into a grid of 1000 equidistant points, and then optimize the prior using a gradient ascent algorithm. See \prettyref{fig:least_fav_prior} for the plot of the least favorable prior.
		
		To emulate the Bayesian oracle, we consider EB estimators based on three methods of estimating the prior $G$ from the data:
	\begin{itemize}
		\item ~Nonparametric maximum likelihood estimator (NPMLE) \cite{KW56}
		
		\item ~Minimum squared Hellinger ($H^2$) distance estimator
		
		\item ~Minimum $\chi^2$-distance estimator.
	\end{itemize}
	These methods are detailed later in \prettyref{sec:EB_results}. We compare their performances with the classical Robbins estimators \cite{Rob51,Rob56} (presented in \eqref{eq:EB_robbins}), the Minimax estimator, and the Naive estimator that directly uses the goals from season 2017-2018 to predict the goals in season 2018-2019 for the same player. The root mean squared error (RMSE) and mean absolute deviation error (MAD) for predicting the hockey goals are presented in \prettyref{tab:minimax_comparison}.\footnote{Given data points $Y_1,\dots,Y_n$ and their predictions $\hat Y_1,\dots,\hat Y_n$ the RMSE is defined as $\sqrt{\frac{1}{n}\sum_{i=1}^n(Y_i-\hat Y_i)^2}$ and the MAD is defined as $\frac{1}{n}\sum_{i=1}^n|Y_i-\hat Y_i|$.} Notably, all EB estimators based on $G$-estimation perform better than the minimax estimator and the Robbins estimator in both error metrics. The above three EB methods based on $G$-estimation also outperform the Naive estimator in terms of RMSE. The performances for the Naive method and our EB methods are similar in the MAD metric. However, note that the construction of the above EB estimators is aimed at estimating the minimum mean squared error estimator. It may be possible to construct EB estimators that mimic the Bayes estimator under absolute error loss, which could provide improved performance guarantees in the MAD metric compared to the Naive method. This is left for future directions.
	
	\begin{table}
		\caption{Minimax vs EB estimators}
		\label{tab:minimax_comparison}
		\begin{tabular*}{\columnwidth}{@{\extracolsep\fill}lllllll@{\extracolsep\fill}}
			\hline
			Methods & Robbins & Minimax & { Naive} & minimum-$H^2$ & NPMLE & minimum-$\chi^2$ \\
			\hline
			RMSE & 15.59 & 8.62 & { 6.19} & 6.02 & 6.04 & 6.05 \\
			MAD & 6.64 & 7.54 & { 4.35} & 4.37 & 4.38 & 4.39 \\
			\hline
		\end{tabular*}
	\end{table}
	
	Getting back to the literature review, there are two main avenues to solving the EB problem:		
	\begin{itemize}
		\item $f$-modeling: Construct an approximate Bayes estimator by approximating the marginal density. For example, the Robbins estimator \cite{Rob56}
		is a plug-in estimate of \prettyref{eq:EB_bayes_est} replacing the true $f_G$ with the empirical distribution, leading to
		\begin{align}\label{eq:EB_robbins}
			\hat\theta_j = \hat\theta_{\text{Robbins}}(Y_j|Y_1,\dots,Y_n)
			\eqdef(Y_j+1){N(Y_j+1)\over N(Y_j)},\quad
			N(y)\eqdef |\{i\in [n]:Y_i=y\}|.
		\end{align} 
		
		\item $g$-modeling: We first obtain an estimate ${\hatGn}$ of the prior $G$ from $Y^n$ and then apply the corresponding Bayes estimator formula $\hat \theta_{{\hatGn}}(Y_j)$. Examples of ${\hatGn}$ include the celebrated NPMLE method mentioned above
		\begin{align}
			{\hatGn}=\argmax_G \frac{1}{n} \sum_{i=1}^n \log f_G(Y_i)
			\label{eq:NPMLE}
		\end{align}
		where the maximization is over all priors on $\reals_+$ ( unconstrained NPMLE). 
		When additional information about the prior is available (e.g., compactly supported), it is convenient to incorporate these constraints into the above optimization, leading to a constrained NPMLE.		 	
	\end{itemize}
	
	In a nutshell, both $f$-modeling and $g$-modeling rely on an estimate of the population density $f_G$; the difference is that the former applies improper density estimate such as the empirical distribution or kernel density estimate (see, e.g., \cite{li2005convergence,brown2009nonparametric,Z09} for Gaussian models), while the latter applies \emph{proper} density estimate of the form $f_{\hat G}$.

In recent years, there have been significant advances in the theoretical analysis of $f$-modeling EB estimators for the Poisson model, specifically, the Robbins method. For compactly supported priors, \cite{BGR13} showed that with Poisson sampling (replacing the sample size $n$ by $\Poi(n)$), the Robbins estimator achieves a $O\left({(\log n)^2\over n (\log\log n)^2}\right)$ regret for estimating each $\theta_i$. Later \cite{PW21} showed that the same bound holds with fixed sample size $n$ and established the optimality of the Robbins estimator by proving a matching minimax lower bound. For the class of subexponential priors, for estimating each $\theta_i$, the Robbins estimator also achieves optimal minimax regret $\Theta\left({(\log n)^3\over n}\right)$.

On the other hand, despite its simplicity and optimality, it has long been recognized that the Robbins method often produces unstable estimates in practice. This occurs particularly for that $ y$ which appears a few times or none whatsoever, so that $N(y)$ is small or zero. Thus, unless $N(y+1)$ is also small, the formula 
\prettyref{eq:EB_robbins}	produces exceptionally large value of 
$\hat\theta_{\text{Robbins}}(y)$. In addition, if $N(y+1)=0$ (e.g., when $y\geq \max\{Y_1,\ldots,Y_n\}$), we have $\hat\theta_{\text{Robbins}}(y)=0$ irrespective of any existing information about $y$, which is at odds with the fact that the Bayes estimator $\hat\theta_{G}(y)$ is always monotonically increasing in $y$ for any $G$ \cite{HS83}.
These issues of the Robbins estimator have been well-documented and discussed in the literature; see, for example, \cite[Section 1]{maritz1968smooth} and \cite[Section 1.9]{maritz2018empirical} for a finite-sample study and \cite[Section 6.1]{efron2021computer} for the destabilized behavior of Robbins estimator in practice (e.g., in analyzing insurance claims data).
To alleviate the shortcomings of the Robbins estimator, a number of modifications have been proposed \cite{maritz1968smooth,BGR13} that enforce smoothness or monotonicity; nevertheless, it is unclear if they still retain the regret optimality of the Robbins method. This raises the question of whether it is possible to construct a well-behaved EB estimator that is provably optimal in terms of regret.

In this paper, we answer this question in the positive. This is accomplished by a class of $g$-modeling EB estimators, which are free from the unstable behavior of the Robbins estimator, thanks to their Bayesian form, which guarantees monotonicity among many other desirable properties. 
The prior is learned using the \emph{minimum-distance} method, including the NPMLE \prettyref{eq:NPMLE} as a special case.		
Introduced in the pioneering works \cite{wolfowitz1953estimation,wolfowitz1954estimation,wolfowitz1957minimum}, 
the minimum-distance method aims to find the best fit \emph{in class} to the data with respect to a given distance.
As such, it is well-suited for estimating the prior, and the obtained density estimate is \emph{proper} and of the desired mixture type.

\begin{figure}[t]
		\centering
		\begin{minipage}{0.5\textwidth}
			\begin{center}
				\includegraphics[ height=6cm]{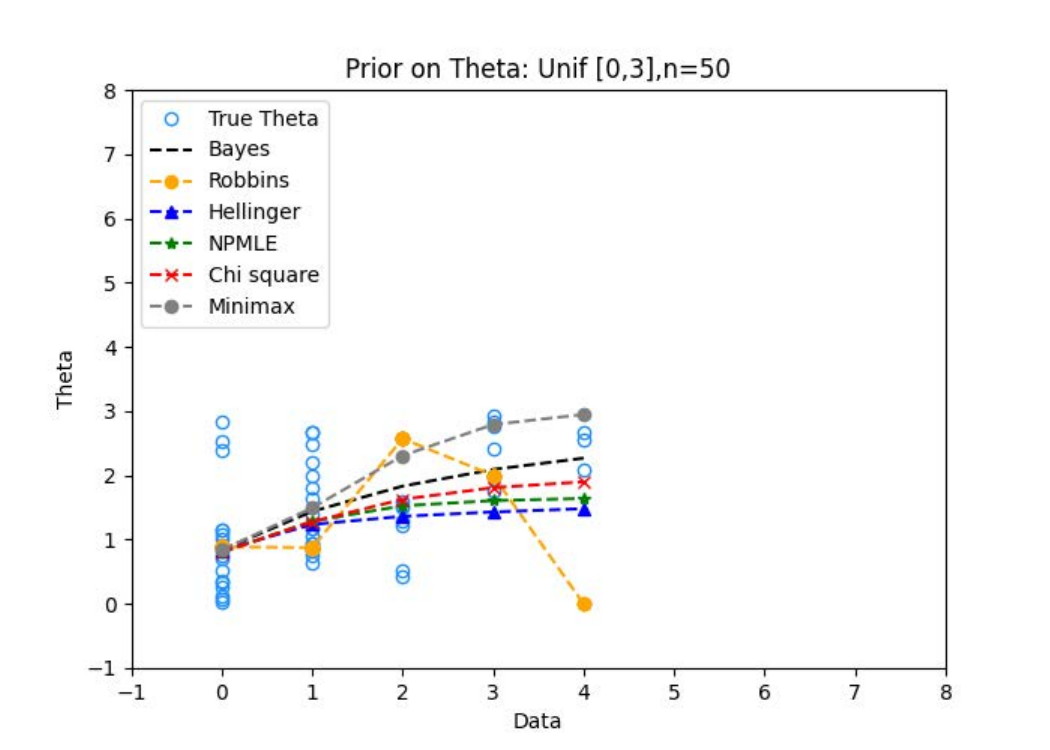}
			\end{center}
		\end{minipage}\hfill
		\begin{minipage}{0.5\textwidth}
			\begin{center}
				\includegraphics[ height=6cm]{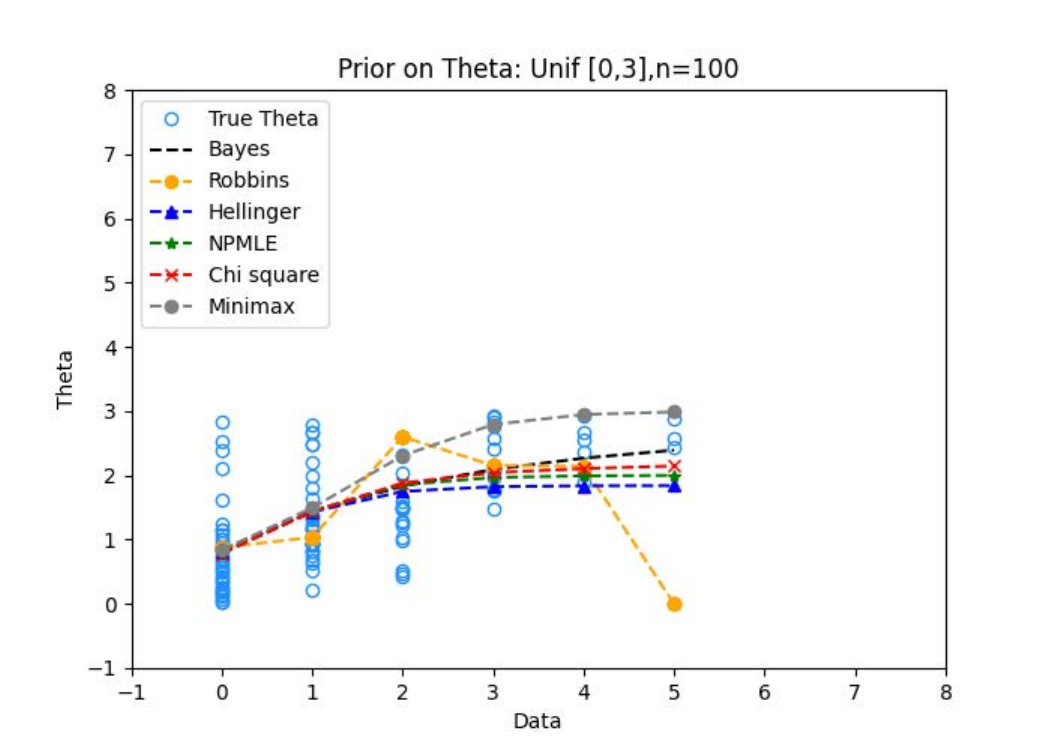}
			\end{center}		
		\end{minipage}
		\begin{minipage}{0.5\textwidth}
			\begin{center}
				\includegraphics[ height=6cm]{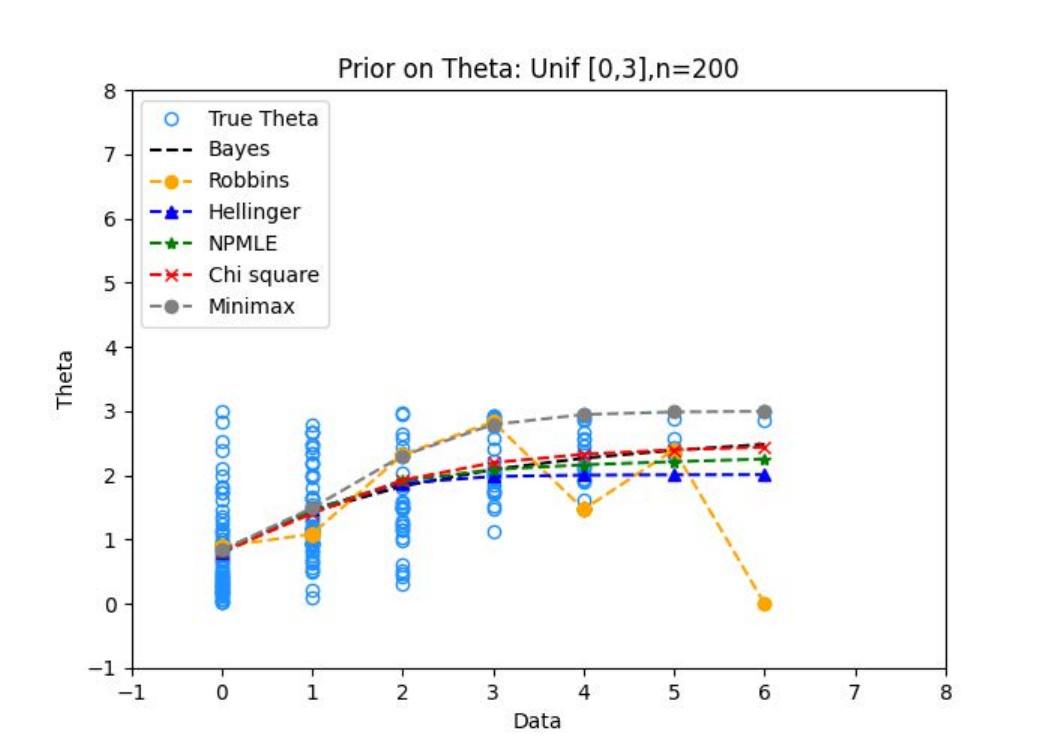}
			\end{center}
		\end{minipage}\hfill
		\begin{minipage}{0.5\textwidth}
			\begin{center}
				\includegraphics[ height=6cm]{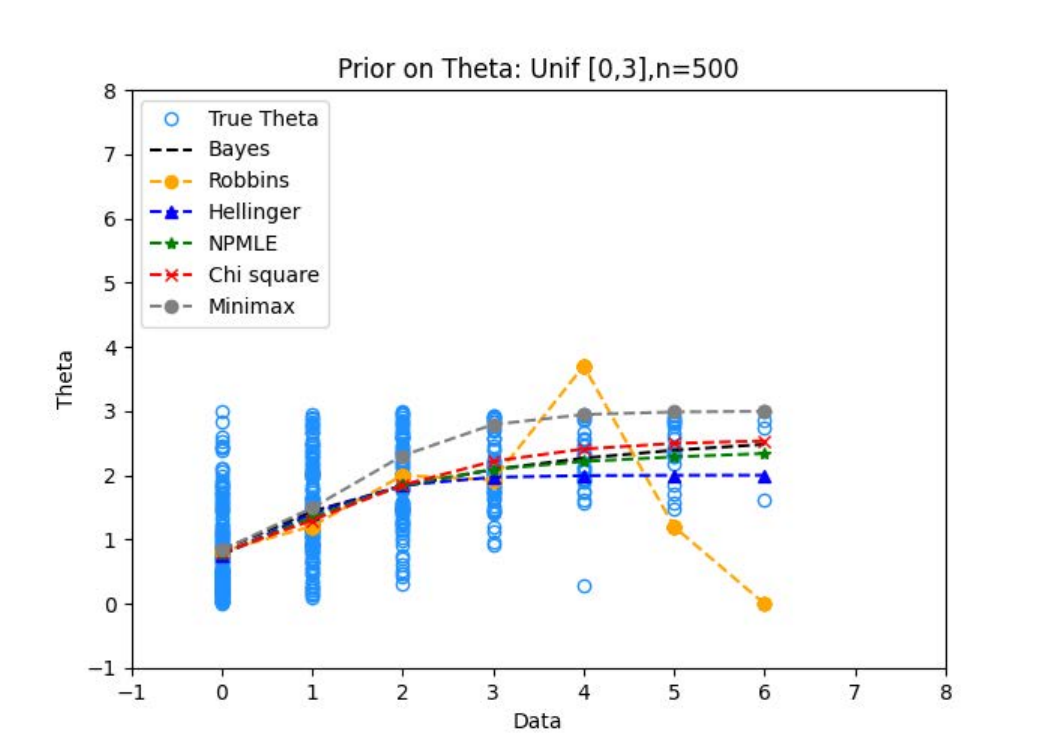}
			\end{center}		
		\end{minipage}
		\caption{Comparison of Robbins estimator with different minimum-distance EB estimators. Here the latent $\theta_i \iiddistr \Unif[0,3]$ and the observation $Y_i\inddistr \Poi(\theta_i)$, for $i=1,\ldots,n$. We plot $\hat \theta(Y_i)$ against $Y_i$ for various EB estimators $\hat \theta$. For reference, we also plot the true value $\theta_i$ and the Bayes estimator $\hat\theta_G(Y_i)$. The sample sizes are $n=50,100,200,500$.}
		\label{fig:Rob-vs-NPMLE}
	\end{figure}

As a concrete example, we consider a simple uniform prior and compare the numerical performance of Robbins and three prototypical examples of minimum-distance estimators of $G$, with respect to the Kullback-Leibler (KL) divergence (i.e., the NPMLE), the Hellinger distance, and the $\chi^2$-divergence, respectively (see \prettyref{sec:min-dist-est} for the formal definitions). As evident in
\prettyref{fig:Rob-vs-NPMLE}, the minimum-distance EB estimators provide a much more consistent approximation of the Bayes estimator compared to the Robbins estimator and the minimax estimator, which can be calculated using the least favorable prior obtained via \eqref{eq:least-favorable-prior-calc} with $h=3$. This advantage is even more pronounced for unbounded priors (cf.~\prettyref{fig:unbounded-priors} in \prettyref{sec:simulation-unbounded}); see also \prettyref{fig:real} for a real-world example where EB methodology is applied to a prediction task with sports data. Notably, in multidimensional settings, such minimum distance based EB methodologies are difficult to implement in practice as they are computationally expensive even in fixed dimensions. However, we propose that the unidimensional EB methodology can also be employed to provide improved analyses in multidimensional setups. To demonstrate the above, we considered a regression problem based on simulations. We demonstrate that the performance of an ordinary least squares (OLS) method can be significantly improved by pre-processing the individual covariate columns using minimum-distance EB filters before supplying them to the algorithm.

	The superior performance of minimum-distance EB estimators in practice is also justified by theory. In addition to characterizing their structural properties (existence, uniqueness, discreteness) in the Poisson model, we show that, under appropriate conditions on the distance functional, their regret is minimax optimal for both compactly supported and subexponential priors. This is accomplished by first proving the optimality of minimum-distance estimate for density estimation in Hellinger distance, then establishing a generic regret upper bound for $g$-modeling EB estimators in terms of the Hellinger error of the corresponding density estimates. We also extend the theoretical analyses to a multidimensional Poisson models.

	\subsection{Related works}
	
	Searching for a stable and smooth alternative to the classical Robbins method for the Poisson EB problem has a long history. \cite{maritz1966smooth} was one of the proponents of using $g$-modeling estimators to resolve this problem. The author considered modeling the prior using the Gamma distribution and estimated the scale and shape parameters using a $\chi^2$-distance minimization. This is a parametric approach as opposed to the nonparametric approach in this paper.
	Based on the monotonicity of the Bayes estimator, 
	\cite{maritz1969empirical} 
	used non-decreasing polynomials to approximate the Bayes oracle (recently, similar isotonic regression based EB approaches have also been used to estimate the Bayes oracle in other models, e.g., see \cite{barbehenn2022nonparametric} for an example in the Gaussian mean estimation problem). Nonetheless, it is unclear whether these regression-based approaches for estimating the Bayes oracle can be used directly to draw any conclusions about estimating the marginal data distribution, which is also one of the primary focuses of our work. \cite{lemon1969empirical} proposed an iterative method of estimating the prior by first using the empirical distribution of the training sample $Y^n$ and then using corresponding posterior means of the $\theta_i$'s to denoise. 
	In a similar vein, \cite{bennett1972continuous} assumed the existence of a density of the prior distribution and used the kernel method to approximate the prior. For a detailed exposition of other smooth EB methods, see \cite{maritz2018empirical}. 
	However, none of these methods has theoretical guarantees in terms of the regret for the nonparametric class of priors considered in the present paper.

	Applying NPMLE to estimate the mixture distribution has been well-studied in the literature. \cite{KW56} was one of the preliminary papers to prove the consistency of the NPMLE, which was subsequently extended in \cite{heckman1984method,jewell1982mixtures,lambert1984asymptotic,pfanzagl1988consistency}; for a more recent discussion, see \cite{chen2017consistency}. In the present paper, we focus on the Poisson mixture model and sharpen these results by obtaining the optimal rate of convergence for the NPMLE.
	In addition to the aforementioned statistical results, structural understanding of the NPMLE (existence, uniqueness, and discreteness) has been obtained in \cite{S76,jewell1982mixtures,L83general,L83,L95} for the general univariate exponential family. We extend these structural results to a class of minimum-distance estimators for Poisson mixture models following \cite{S76}. Finally, we mention the recent work \cite{miao2021fisher}, which explored the application of NPMLE in a related scenario of heterogeneous Poisson mixtures.

	Initial work on applying NPMLE for EB estimation was carried out in \cite{laird1982empirical} for the Binomial and the normal location models, and the analysis is primarily numerical. For theoretical results, \cite{ghosal.vdv,Z09} analyzed the Hellinger risk of NPMLE-based mixture density estimates, which forms the basis of the analysis of NPMLE for EB estimation in \cite{JZ09}. The resulting regret bounds, though state-of-the-art, still differ from the minimax lower bounds in \cite{PW21} by logarithmic factors for both the classes of compactly supported and subgaussian priors. This is because (a) the density estimation analysis in \cite{Z09} is potentially suboptimal compared to the lower bounds in \cite{K14};
	(b) the Fourier-analytic reduction from the Hellinger distance for mixture density to regret in \cite{JZ09} is loose.
	In comparison, density estimation and regret bounds are optimal in this paper with exact logarithmic factors. This can be attributed to the discrete nature of the Poisson model, which allows a simple truncation-based analysis to suffice for light-tailed priors. These sharp results are also generalized from the NPMLE-based EB estimator to the minimum-distance estimators.
	
	The minimum distance based analysis for the Poisson model we consider here can be easily extended to other members of the exponential family of distribution. In a recent follow-up work \cite{jana2023empirical}, it was shown that the Bayes estimator for various distributions, including the Geometric distribution, Negative Binomial distribution, and Exponential distribution, can be represented similarly in terms of simple functions of the prior distribution. Our minimum distance methodology can then be extended to estimate the prior distribution, and similar regret analysis can be generalized to these discrete models.

	\subsection{Organization}
	The rest of the paper is organized as follows. In \prettyref{sec:EB_results}, we introduce the class of minimum distance estimators and identify conditions on the distance function that guarantee the minimizer's existence and uniqueness. The theoretical guarantees in density estimation and regret are presented in \prettyref{thm:EB_density-H2} and \prettyref{thm:EB_main} therein. The proof sketches of these results are presented in \prettyref{sec:EB_density-estimation} and \prettyref{sec:EB_regret-upper-bound}, respectively. In \prettyref{sec:numerical-expt}, we present an algorithm for computing minimum-distance estimators in the one-dimensional setting. We also study their numerical performance in empirical Bayes estimation with simulated and real datasets. In \prettyref{sec:results-multidim}, we mention our theoretical results in a multi-dimensional setting. For all the other related details of proofs, see the Appendix.
	\subsection{Notations}
	\label{sec:notation}
	Denote by $\integers_+$ (resp.~$\reals_+$) the set of non-negative integers (resp.~real numbers). 
	For a Borel measurable subset $\Theta\subset\reals$, let $\calP(\Theta)$ be the collection of all probability measures on $\Theta$. 
	For any $\theta\in\reals_+$ let $\delta_\theta$ denote the Dirac measure at $\theta$. 
	Denote by $\subexpo(s)$ the set of all $s$-subexponential distributions on $\reals_+$: $\subexpo(s)=\sth{G:G([t,\infty)])\leq 2e^{-t/s},\forall t>0}$. 
	Let $Y_i\sim \Poi(\theta_i)$ for $i=1,\dots,n$ and $Y\sim\Poi(\theta)$, with $\theta_1,\dots,\theta_n,\theta\iiddistr G$. This also implies  $Y_1,\dots,Y_n,Y\iiddistr f_G$ where $f_G$ is the mixture distribution defined in \eqref{eq:EB_mixture}.
	Let $\EE_G$ and $\PP_G$ denote the expectation and probability where the true mixing distribution is $G$. Define $Y_{\max}=\max_{i\in [n]}Y_i, Y_{\min}=\min_{i\in [n]}Y_i$.

	\section{Problem formulation and results}
	\label{sec:EB_results}
	\subsection{Minimum-distance estimators}\label{sec:min-dist-est}
	Denote by $\calP(\integers_+)$ the collection of probability distributions (pmfs) on $\mathbb Z_+$. 
	We call $\dist: \calP(\integers_+)\times\calP(\integers_+) \to\reals_+$ a \emph{generalized distance} if
	$\dist(p\|q)\geq 0$ for any $p,q\in\calP(\integers_+)$, with equality if and only if $p=q$. 
	Note that any metric or $f$-divergence \cite{Csiszar67} qualifies as a generalized distance.

	The minimum-distance\footnote{We adopt this conventional terminology even when $\dist$ need not be a distance.} methodology aims to find the closest fit in the model class to the data. While it is widely used and well-studied in parametric models \cite{beran1977minimum,berkson1955maximum,pollard1980minimum,bolthausen1977convergence,millar1984general}, it is also useful in nonparametric settings such as mixture models.
	Denote by
	\begin{align}\label{eq:EB_pmf-empirical}
		p^\sfem_n=\frac 1n\sum_{i=1}^n\delta_{Y_i}
	\end{align} 
	the empirical distribution of the sample $Y_1,\dots,Y_n$. 
	The minimum-distance estimator for the mixing distribution with respect to $\dist$, over some target class of distributions $\calG$, is
	\begin{equation}
		\hatGn \in \argmin_{Q\in\calG}\dist(p^\sfem_n\|f_Q).
		\label{eq:mindist}
	\end{equation}
	Note that in our analysis we also allow $\calG = \calP(\reals_+)$, the set of all probability distributions on the nonnegative real line. In such cases, for ease of notations, we will replace $\argmin_{Q\in\calP(\reals_+)}\dist(p^\sfem_n\|f_Q)$ with $\argmin_Q\dist(p^\sfem_n\|f_Q)$. Primary examples of minimum-distance estimators considered in this paper include the following  
	\begin{itemize}
		\item ~Maximum likelihood: $\dist(p\|q) = \mathrm{KL}(p\|q)\triangleq \sum_{y \geq 0} p(y)\log \frac{p(y)}{q(y)}$ is the KL divergence. 
		In this case, one can verify that the minimum-KL estimator coincides with the NPMLE \prettyref{eq:NPMLE}.
		
		\item ~Minimum-Hellinger estimator: $\dist(p\|q) = H^2(p,q)\triangleq\sum_{y\geq 0}\pth{\sqrt{p(y)}-\sqrt{q(y)}}^2$ is the squared Hellinger distance.
		
		\item ~Minimum-$\chi^2$ estimator: $\dist(p\|q) = \chi^2(p\|q)\triangleq\sum_{y \geq 0}{(p(y)-q(y))^2\over q(y)}$ is the $\chi^2$-divergence.
	\end{itemize}
	Note that there are other minimum-distance estimators previously studied for the Gaussian mixture model, e.g., those respect to $L_p$-distance of the CDFs,
	aiming at estimating the mixing distribution \cite{DK68,Chen95,HK2015,edelman1988estimation}. These are outside the scope of the theory developed here. 
	
	In general, the solution to \prettyref{eq:mindist} need not be unique; nevertheless, for the Poisson mixture model, the uniqueness is guaranteed provided that 
	the generalized distance $\dist$ admits the following decomposition:
	
	\begin{assumption}
		\label{pt:EB_t-ell} 
		There exist maps $t:\calP(\integers_+)\to \reals$ and $ \ell:\reals^2\to\reals$ such that for any two distributions $q_1,q_2\in\calP(\integers_+)$
		$$\dist\pth{q_1\| q_2}= t(q_1)+\sum_{y\geq 0}\ell(q_1(y),q_2(y)),$$
		where $b \mapsto \ell(a,b)$ is strictly decreasing and strictly convex for $a> 0$ and $\ell(0,b)=0$ for $b\geq 0$. 
	\end{assumption}

	The following theorem guarantees the existence, uniqueness, and discreteness of both unconstrained and support-constrained minimum-distance estimators.
	For the special case of unconstrained NPMLE this result was previously shown by \cite{S76} and later extended to all one-dimensional exponential family \cite{L95}.
	\begin{theorem} \label{thm:existence-uniqueness}
		Let $\dist$ satisfy \prettyref{pt:EB_t-ell}. Let $p$ be a probability distribution on $\integers_+$ with support size $m$. Then for any $h>0$, the constrained solution $\argmin_{Q\in\calP([0,h])} \dist(p\|f_{Q})$ exist uniquely and is a discrete distribution with support size at most $m$.
		Furthermore, the same conclusion also applies to the unconstrained solution $\argmin_{Q\in\calP(\reals_+)} \dist(p\|f_{Q})$, which in addition is supported on $[\min_{i=1,\ldots,m} y_i,\max_{i=1,\ldots,m} y_i]$, where $\{y_1,\dots,y_m\}$ is the support of $p$.
	\end{theorem}

	To analyze the statistical performance of minimum-distance estimators, we impose the following regulatory condition on the generalized distance $\dist$:
	\begin{assumption}
		\label{pt:EB_sandwich}
		There exist absolute constants $c_1,c_2>0$ such that for pmfs $q_1,q_2$ on $\integers_+$
		\begin{align}\label{eq:EB_Assupmtion_2}
			c_1H^2(q_1,q_2)\leq \dist(q_1\|q_2)\leq c_2\chi^2(q_1\|q_2)
		\end{align}
	\end{assumption}
	Major examples of generalized distance satisfying  Assumptions \ref{pt:EB_t-ell} and \ref{pt:EB_sandwich} include the KL divergence, squared Hellinger distance, and $\chi^2$-divergence. This follows from noting that $2H^2\leq\text{KL}\leq\chi^2$ and each of them satisfies the decomposition \prettyref{pt:EB_t-ell}: for squared Hellinger $t\equiv 2, \ell(a,b)=-2\sqrt{ab}$, for KL divergence $t\equiv 0, \ell(a,b)=a\log\frac ab$, for $\chi^2$-divergence $t\equiv -1, \ell(a,b)={a^2\over b}$. 
	On the other hand, total variation (TV) satisfies neither \prettyref{pt:EB_t-ell} nor \ref{pt:EB_sandwich} so the theory in the present paper does not apply to the minimum-TV estimator.
	
	\begin{remark}
		Before proceeding further, note the following argument in the context of \prettyref{thm:existence-uniqueness} (which is a deterministic result) to exclude $G=\delta_0$ as a possible choice in the analysis, and the situation where all $y_i$-s are zero. The remark applies for the rest of the paper to exclude $G=\delta_0$ as a possible choice in the analysis. Whenever we need to divide with $f_G(y)$, the choice $G=\delta_0$ will lead to technical difficulties as $f_G$ will be degenerate at zero as well. Such division occurs, for example, in the proof of \prettyref{thm:EB_density-H2}, where we use a bound based on the $\chi^2$-divergence between the empirical distribution and $f_G$ to control the behavior of the general minimum distance estimator. Consider the following cases:
			\begin{itemize}
				\item ~$y_1=...=y_n=0$. Then clearly $\hat G=\delta_0$ is the unique NPMLE solution. 
			\item ~$y_i>0$ for some $i\in [n]$. Then clearly $\delta_0$ is not the data generating distribution. As a result, we may assume $f_G$ to be fully supported on $\integers_+$ for priors $G$ in the analysis.
			\end{itemize}
	\end{remark}

	\subsection{Main results}
	
	In this section we state the statistical guarantee	for the minimum-distance estimator $\hat G$ defined in the previous section. Our main results are two-fold (both minimax optimal): 
	\begin{enumerate}
		\item ~Density estimation, in terms of the Hellinger distance $f_{\hat G}$ and the true mixture $f_{G}$;
		\item ~Empirical Bayes, in terms of the regret of the Bayes estimator with the learned prior $\hat G$.
	\end{enumerate}
	As mentioned in \prettyref{sec:EB_introduction}, the regret analysis in fact relies on bounding the density estimation error. 
	We start with the result for density estimation. 
	Recall from \prettyref{sec:notation} $\calP([0,h])$ and $\subexpo(s)$ denote the class of compactly supported and subexponential priors respectively.
	
	\begin{theorem}[{Density estimation}]\label{thm:EB_density-H2}
		Let $\dist$ satisfy Assumption~\ref{pt:EB_t-ell} and Assumption~\ref{pt:EB_sandwich}. Suppose that $\hat G$ is the unconstrained minimum distance estimator
		\begin{align}
			\label{eq:est-unconstrained}
		\hat G = \argmin_{Q\in \calP(\reals^+)} \dist(p_n^\sfem\|f_{Q}).
		\end{align}
		Then there exist constants $c_1,c_2$ such that the following holds.
		\begin{enumerate}[label=(\alph*)~]
			\item  $\sup_{G\in \calP[0,h]}\EE\qth{H^2(f_{G},f_{{\hatGn}})}\leq {c_1\over n}\cdot \min\sth{\max\{1,h\}{\log n\over \log\log n},h+\log n}$ for any $n\geq 3$.
			\item $\sup_{G\in\subexpo(s)}\EE\qth{H^2(f_{G},f_{{\hatGn}})}\leq 
			{c_2\max\{1,s\}\over n}{\log n}$ for any $n\geq 2$.
		\end{enumerate}

	\end{theorem}
	\begin{remark}\label{rmk:density}
		\begin{enumerate}[label=(\roman*)~]
			
			\item 
				If the prior $G$ is allowed to be any distribution on $\reals^+$, then neither density estimation nor empirical Bayes estimation is possible. 
				This fact is well-known for the Gaussian mixture model \cite{Z09,suresh2014near}.
				To see this, fix a quadratically spaced grid 
				$\{\theta_1,\ldots,\theta_n\}$ where $\theta_i = i^2 \cdot (\log n)^{10}$.
				Consider a prior $G = \frac{1}{n} \sum_{i=1}^n \text{Uniform}(\theta_i,\theta_i+\alpha_i)$, where $\alpha_i$'s are iid drawn from $\text{Uniform}(0,1)$. 
				In other words, the prior $G$ is a uniform mixture over $n$ clusters each of which has an $O(1)$ spread. The quadratic grid is chosen so that the spacing $\theta_{i+1}-\theta_i$ far exceeds the standard deviation $\sqrt{\theta_i}$ so that with high probability we know which cluster each data point $Y_i$ is drawn from. 
				However, there is no enough information to estimate the parameters $\alpha_i$'s consistently because on average we only observe one sample for each cluster. 				
				As such, in order to obtain uniform error bound as in \prettyref{thm:EB_density-H2} (and later in \prettyref{thm:EB_main} for regret), it is necessary to restrict the priors to a subclass.

			\item It has been shown recently in \cite[Theorem 21]{PW21} that for any constant $h,s$, the minimax squared Hellinger density estimation errors are at least $\Omega\left(\log n\over n\log\log n\right)$ and $\Omega\left(\log n\over n\right)$ for priors in the class $\calP([0,h])$ and $\subexpo(s)$, respectively. This establishes the minimax optimality of our minimum-distance density estimates.

				\item In the shape-constrained density estimation literature, e.g., \cite{koenker2018shape}, there has also been some interest in replacing the maximum likelihood fitting criteria with other divergences. In this context, the motivation was the desire to impose weaker concavity constraints than log-concavity while still preserving the underlying convexity of the variational formulation of the problem. Thus, in that setting, divergences were dictated by the form of the concavity constraints. Such convexity assumption is often a necessary criterion to guarantee the validity of divergences; for example, in the case of general Bregman divergences \cite{jana2019characterization,ray2022characterizing}, the functional density power divergences \cite{ray2022characterizing}, and the convexity of the loss function might result in practical benefits, such as the efficiency of minimum divergence estimators \cite[Section 2]{lindsay1994efficiency},\cite{basu2011statistical}. In contrast, our work imposes such convexity constraints on the divergences mainly to guarantee the uniqueness and existence of the estimators. Our theoretical analysis does not require convexity constraints as we primarily aim to establish the finite sample error rate of our minimum distance estimators. Using different convexity structures may result in different multiplicative constants in the minimax error rates, which is beyond the scope of our current work. A detailed study of the convexity structures might help to differentiate between the performances of the minimum-distance estimators, and we leave it for future directions.
				\item Our current theoretical results become vacuous when $h$ approaches a similar order of magnitude to $n$. However, this is mostly due to technicalities in the proof of our upper bound where we aimed to perform a uniform analysis for all priors over $[0,h]$ for a constant $h$. We can revise our analysis to achieve consistency guarantees of the estimators for larger values of $h$, however the general treatment is beyond the scope of the current work. As an example, we detail below the performance guarantee of the unconstrained minimum distance estimator $\hat G$ (the estimator $\hat G$ does not know $G$ or $h$) when the data-generating prior $G$ is degenerate at some $h$, as suggested by one of the reviewers. Our analysis includes the the case when $h$ can be larger than $n$, but significantly smaller than $n^2/\log n$. Note that the degenerate prior corresponds to the data distribution $Y_1\dots,Y_n\simiid \Poi(h)$. We will show that the minimum distance-based estimator $f_{\hat G}$ of $f_G$ should consistently estimate $\Poi(h)$. Revisiting the proof of Theorem 2 and using that the unconstrained minimum distance estimator $\hat G$ is supported on $[Y_{\min},Y_{\max}]$, we get
				\begin{align*}
					H^2(f_G,f_{\hat G})
					&\leq 
					2\pth{H^2(p_n^\sfem,f_{\hat G})+ H^2(p_n^\sfem,f_{G})}
					\nonumber\\
					&\leq
					\frac 2{c_1}\pth{\dist(p_n^\sfem\|f_{\hat G})+\dist(p_n^\sfem\|f_{G})}
					\leq \frac 4{c_1} \dist(p_n^\sfem\|f_{G})
					\leq
					{4c_2\over c_1}\chi^2(p_n^\sfem\|f_{G}),
				\end{align*}
				and hence, for $K_1,K_2$ to be chosen later, we have
				\begin{align}\label{eq:r1}
					\EE\qth{H^2(f_G,f_{\hat G})} 
					&=\EE\qth{H^2(f_G,f_{\hat G})\indc{Y_{\min}<K_1 \text{ or } Y_{\max}>K_2}}
					+\EE\qth{H^2(f_G,f_{\hat G})\indc{K_1<Y_{\min}\leq Y_{\max}<K_2}}
					\nonumber\\
					&\leq 
					4\PP\qth{Y_{\min}<K_1 \text{ or } Y_{\max}>K_2}
					+{4c_2\over c_1}\EE\qth{\chi^2(p_n^\sfem\|f_G)\indc{K_1<Y_{\min}\leq Y_{\max}<K_2}}
					\nonumber\\
					&\leq 4n\PP\qth{Y<K_1}+4n\PP\qth{Y>K_2}
					+{4c_2\over c_1}\EE\qth{\chi^2(p_n^\sfem\|f_G)\indc{K_1<Y_{\min}\leq Y_{\max}<K_2}}
				\end{align}
				for a random variable $Y\sim \Poi(h)$, where the last inequality used union bounds. We can bound the rightmost term in the above expression as
				\begin{align}
					&\EE\qth{\chi^2(p_n^\sfem\|f_G)\indc{K_1<Y_{\min}\leq Y_{\max}<K_2}}
					\nonumber\\
					&=\sum_{y\in [K_1,K_2]}{\EE\qth{(p_n^\sfem(y)-f_G(y))^2}\over f_G(y)}
					+\sum_{y<K_1 \text{ or } y>K_2} f_G(y)\PP[K_1<Y_{\min}\leq Y_{\max}<K_2]
					\nonumber\\
					&= 
					\sum_{y\in [K_1,K_2]}{1-f_G(y)\over n}
					+\PP[K_1<Y_{\min}\leq Y_{\max}<K_2] 
					\cdot \PP_{Y\sim \Poi(h)}\qth{Y<K_1 \text{ or } Y>K_2}
					\nonumber\\
					&\leq
					{K_2-K_1\over n} + \PP_{Y\sim \Poi(h)}\qth{Y<K_1}+ \PP_{Y\sim \Poi(h)}\qth{Y>K_2}.
				\end{align}
				In view of \eqref{eq:r1} the last display implies 
				\begin{align}\label{eq:r2}
					\EE\qth{H^2(f_G,f_{\hat G})}
					\leq {4c_2(K_2-K_1)\over c_1n}+(4n+{4c_2\over c_1})(\PP\qth{Y<K_1}+\PP\qth{Y>K_2})
				\end{align}
					We choose
					$$
					K_1 = \max\{0, h-2\sqrt{h\log n}\},\quad K_2=h+3\sqrt{h\log n}.
					$$
					Then we can use
					\prettyref{lmm:anru-zhang-results-poisson} to bound $\PP\qth{Y<K_1},\PP[Y>K_2]$ as long as $h\geq c\log n$ for a large enough $c$. The application gives us $\PP\qth{Y<K_1},\PP[Y>K_2]\leq \frac 2{n^2}$.
					Plugging the above choice in \eqref{eq:r2}, we get for a constant $c>0$
					$$
					\EE\qth{H^2(f_G,f_{\hat G})}
					\leq {c\sqrt{h\log n}\over n}.
					$$
					This implies that the minimum distance estimator is consistent given $h$ is significantly smaller than $n^2/\log n$.
		\end{enumerate}
	\end{remark}

	Next, we turn to the problem of estimating $\theta_1,\ldots,\theta_n$ from $Y_1,\ldots,Y_n$, under the squared error loss, using the empirical Bayes methodology. In this work, we study the estimation guarantees of the $g$-modeling type estimators. Notably, to produce an estimator $\hat\theta_j$ of $\theta_j$, we use the observations $Y^{-j} = (Y_1,\dots,Y_{j-1},Y_{j+1},\dots,Y_n)$ to approximate $G$ and then plug it in the formula of the Bayes estimator $\hat\theta_{\hatGn}(Y_j)$ in \eqref{eq:EB_bayes_est}. Given any class of distributions $\calG$ and any distribution estimator strategy characterized by $\hat G$, define the total regret as its worst-case excess risk over the Bayes error:
\begin{align}
		{\sf TotRegret}_n(\hat G;\calG)&\eqdef \sup_{G\in \calG}\sth{\EE_{G}\qth{\|\hat\theta^n(Y^n)-\theta^n\|^2-n\cdot \mmse(G)}},\\
		\hat\theta_j &= \hat \theta_{\hat G (Y^{-j})}(Y_j),\quad j=1,\dots,n.
\end{align}
where $\mmse(G)$ denotes the minimum mean squared error of estimating $\theta\sim G$ based on a single observation $Y\sim f_\theta$, i.e., the Bayes risk
\begin{align}
	\mmse(G)\eqdef \inf_{\hat \theta}~\EE_G
	\qth{\pth{\hat \theta(Y)-\theta}^2}
	=\EE_G\qth{(\hat\theta_G(Y)-\theta)^2}.
\end{align} 
In addition, define the problem of quantifying the individual regret for the estimator $\hat G$
\begin{align}
	\Regret_n(\hat G;\calG) \eqdef \sup_{G\in \calG}\sth{\EE_G\qth{(\hat\theta_n(Y^n)-\theta_n)^2}-\mmse(G)},\quad 
	\hat\theta_n(Y^n) = \hat \theta_{\hat G (Y^{n-1})}(Y_n).
\end{align}
Here $Y_1,\dots,Y_{n-1}$ can be viewed as training data which is used to learn the estimator, and then we apply it on a fresh (unseen) data point $Y_n$ to predict $\theta_n$. Turning to the loss function under consideration, it is not difficult to show that the total regret with sample size $n$ can be bounded from above using $n$ times the individual regret with training sample size $n-1$
\begin{align}
	{\sf TotRegret}_n(\hat G;\calG)\leq n\cdot \Regret_{n}(\hat G;\calG).
\end{align}
In view of the above inequality of the total and individual regret functions, we limit ourselves to studying individual regret only, as this will suffice to achieve the desired optimal rates.

Now, we are in a position to describe the main results for empirical Bayes estimation. For an ease of notation, suppose that given a fresh sample $Y\sim \Poi(\theta)$, where $\theta$ is generated from an unknown prior distribution $G$, we want to predict the value of $\theta$ in the squared error loss and training sample to construct the estimator $\hat G$ is given by $Y_1,\dots,Y_n$. Given any estimator ${\hatGn}$ of $G$ we define the regret of the empirical Bayes estimate $\hat \theta_{\hatGn}$ as
	\begin{align}
		\Regret({\hatGn};G)&
		={\EE_G\qth{\pth{\hat\theta_{{\hatGn}}(Y)-\theta}^2}}-\mmse(G)
		\nonumber\\
		&\stepa{=}{\EE_G\qth{\pth{\hat\theta_{{\hatGn}}(Y)-\hat\theta_{ G}(Y)}^2}} 
		= \EE_G\qth{\sum_{y\geq 0}
			(\hat\theta_{{\hatGn}}(y)-\hat\theta_{ G}(y))^2 f_G(y)},
		\label{eq:regret}
	\end{align}
	where the identity (a) followed by using the orthogonality principle: the average risk of any estimator $\hat \theta$ can be decomposed as 
	\begin{align}\label{eq:EB_ortho-relation}
		\EE_G[(\hat\theta-\theta)^2]
		=\mmse(G)+\EE_G[(\hat\theta-\hat\theta_{G})^2].
	\end{align}
	Similarly we define the maximum regret of ${\hatGn}$ over the class of data generating distributions $\calG$ 
	\begin{align}
		\label{eq:reg-worstcase}
		\Regret({\hatGn};\calG)=\sup_{G\in\calG}\Regret({\hatGn};G).
	\end{align}
	Then we have the following estimation guarantees.
	\begin{theorem}[Empirical Bayes]\label{thm:EB_main}
		Let $\dist$ satisfy Assumption~\ref{pt:EB_t-ell} and Assumption~\ref{pt:EB_sandwich}. Suppose that $\hat G$ is the unconstrained minimum distance estimator given in \eqref{eq:est-unconstrained}. 
		Then there exist constants $c_1,c_2$ such that the following holds
		
		\begin{enumerate}[label=(\alph*)~]
			\item  $
			\Regret({\hatGn};\calP([0,h]))
			\leq {c_1\over n}\cdot \min\sth{\max\{1,h\}{\log n\over \log\log n},h+\log n}^3$ for any $n\geq 3$
		\item $\Regret({\hatGn};\subexpo(s))
		\leq \frac {c_2 \cdot \max\{1,s^3\}}n\pth{\log n}^3$
		\end{enumerate}
		In addition, if the data generating distribution $G$ is supported on $[0,h]$ for a constant $h>0$, then the constrained minimum distance estimator $\tilde G$ with access to $h$ achieves an improved risk guarantee
		\begin{align*}
			\tilde G=\argmin_{Q\in \calP[0, h]}\dist(p_n^\sfem\|f_{Q}),
			\
			\Regret({\hatGn};\calP([0,h]))
			\leq {c_1\over n}\cdot \min\sth{\max\{1,h\}{\log n\over \log\log n},h+\log n}^2.
		\end{align*}
	\end{theorem}
	
	\begin{remark}
		\begin{enumerate}[label=(\roman*)~]
			
			\item As mentioned in \prettyref{sec:EB_introduction}, 
			for fixed $h$ and $s$, the above establishes the regret optimality of minimum-distance EB estimators by matching the minimax lower bounds recently shown in \cite[Theorem 1]{PW21}. This minimax optimality was only known for the $f$-modeling-based Robbins estimator.

			\item  When $\dist$ is the KL divergence, the minimum-distance estimator ${\hatGn}=\argmin_Q \mathsf{KL}(p\|q)$ is the NPMLE. This follows from the expansion
			\begin{align*}
				\text{KL}(p^\sfem_n\|f_Q)
				=\sum_{y\geq 0}p^\sfem_n(y)\log{p^\sfem_n(y)\over f_Q(y)}
				=\sum_{y\geq 0}p^\sfem_n(y)\log{p^\sfem_n(y)}
				-\frac 1n\sum_{i=1}^n\log f_Q(Y_i).
			\end{align*}

			\item \prettyref{thm:EB_main} holds for approximate solutions. Consider the following approximate minimum-distance estimators $\hat G$, over some target class of distributions $\calG$, that satisfies
			\begin{align}
				\label{eq:EB_delta-mindist}
				\dist(p^\sfem_n\|f_{{\hatGn}}) \leq 
				\inf_{Q\in\calG} \dist(p^\sfem_n\|f_Q)+\delta.
			\end{align}
			for some $\delta>0$. Then 
			the regret bound for the bounded prior case (resp. subExponential data generating prior case)			
			continues to hold if $\delta\lesssim {\log n\over n\log\log n}$ (resp.~${\log n\over n}$). Note that ${\hatGn}$ is the NPMLE over $\calG$ if $\delta=0$ and $\dist$ is given by KL divergence. In case of NPMLE, \eqref{eq:EB_delta-mindist} translates to an approximate likelihood maximizer ${\hatGn}$ such that
			\begin{align*}
				\frac 1n\sum_{i=1}^n \log f_{{\hatGn}}(Y_i)\geq \argmax_{G\in\calG}\frac 1n\sum_{i=1}^n\log f_G(Y_i)-\delta.
			\end{align*}
			This type of results is well-known in the literature, see, for example, \cite{JZ09,Z09} for the normal location-mixture model.
			\item
			
			Similar to \prettyref{rmk:density}, our results become trivial when $h$ is of similar order to $n$, and this is primarily due to our proof strategy, which attempts to provide a uniform analysis for all priors over $[0,h]$ for a constant $h$. The estimators in practice perform well in many scenarios when $h$ is large. To demonstrate the above, as suggested by one of the reviewers, we consider the prior distribution degenerate at $h=50$ and performed a simulation study with $n=50$ samples. We present below the performance of the unconstrained minimum distance estimators, in terms of the average RMSE and the average MAD metric, out of 25 repetitions. Notably, out estimators outperform both the Robbins estimator and the naive estimator, that uses the data point $Y_i$ to estimate $\theta_i$-s (which are now identically 50).
			
			\begin{table}[h]
				\caption{Performance of the unconstrained estimators with $h=n=50$}
				
				\label{tab:comparison-50}
				\begin{tabular*}{\columnwidth}{@{\extracolsep\fill}llllll@{\extracolsep\fill}}
					\hline
					Methods & Robbins & { Naive} & minimum-$H^2$ & NPMLE & minimum-$\chi^2$ \\
					\hline
					Average RMSE & 43.94 & {6.94} & 0.84 & 1.25 & 2.14 \\
					Average MAD & 32.91 & {5.56} & 0.83 & 1.05 &	1.72
					 \\
					\hline
				\end{tabular*}
			\end{table}

		\end{enumerate}
	\end{remark}

	\section{Proof for density estimation}
	\label{sec:EB_density-estimation}

	The proof of \prettyref{thm:EB_density-H2} is based on a simple truncation idea. 
	It is straightforward to show that the density estimation error for any minimum distance estimator can be bounded from above, within a constant factor, by the expected squared Hellinger distance between the empirical distribution $p^\sfem_n$ and the data-generating distribution $f_G$, which is further bounded by the expected $\chi^2$-distance. 
	The major contribution to $\chi^2(p^\sfem_n\|f_G)$ comes from the ``effective support" of $f_G$, outside of which the total probability is $o(\frac 1n)$. For the the prior classes $\calP([0,h])$ and $\subexpo(s)$, the Poisson mixture $f_G$ is effectively supported on $\{0,\ldots,O({\log n\over \log\log n})\}$ and $\{0,\ldots,O(\log n)\}$. Each point in the effective support contributes $\frac 1n$ to $\chi^2(p^\sfem_n\|f_G)$ from which our results follow.

	\begin{proof}[Proof of \prettyref{thm:EB_density-H2}]
		
		For any integer $K\geq 1$ and distribution $G$ denote 
		\begin{align}\label{eq:ep-K(G)}
			\epsilon_K(G)\eqdef\PP\qth{Y\geq K}=\sum_{y=K}^\infty f_G(y)
		\end{align}
		Note that $\dist$ satisfies \prettyref{pt:EB_sandwich}, namely \eqref{eq:EB_Assupmtion_2}. 
		We first prove the following general inequality
		\begin{align}\label{eq:EB_general-H2-bound}
			\EE\qth{H^2(f_{G},f_{{\hatGn}})}\leq {4c_2\over c_1}{K\over n}+\pth{{4c_2\over c_1}+2n}\epsilon_K(G).
		\end{align}
		Using the inequality $(a+b)^2\leq 2(a^2+b^2)$, and as ${\hatGn}$ is the minimizer we get
		\begin{align}
			H^2(f_{G},f_{{\hatGn}})
			&\leq \pth{H(p^\sfem_n,f_{{\hatGn}})+H(p^\sfem_n,f_G)}^2
			\leq 2[H^2(p^\sfem_n,f_{{\hatGn}})+H^2(p^\sfem_n,f_{G})]\nonumber\\
			&\leq \frac 2{c_1}{(\dist(p^\sfem_n\|f_{{\hatGn}})+\dist(p^\sfem_n\|f_{G}))}
			\leq \frac 4{c_1}{\dist(p^\sfem_n\|f_{G})}. \label{eq:EB_m6}
		\end{align}
		Define $Y_{\max}
			\eqdef \max_{i=1}^nY_i.$
		as before. 
		Then, bounding ${1\over c_2}d$ by $\chi^2$ we get the following chain
		\begin{align*}
			\frac{1}{c_2} \EE\qth{\dist(p^\sfem_n\|f_{G})\indc{Y_{\max}< K}}
			&\leq \EE\qth{\chi^2(p^\sfem_n\|f_{G})\indc{Y_{\max}< K}} 
			=
			\sum_{y\geq 0}{\EE\qth{(p^\sfem_n(y)-f_{G}(y))^2\indc{Y_{\max}< K}}\over f_{G}(y)}
			\nonumber\\
			&\stepa{=} 
			\sum_{y < K}{\EE\qth{(p^\sfem_n(y)-f_{G}(y))^2\indc{Y_{\max}< K}}\over
				f_{G}(y)}  +\sum_{y\geq K}f_{G}(y) \PP[Y_{\max}<K],
		\end{align*}
		where the last equality follows from the fact that under $\{Y_{\max}<K\}$ we have
		$p^{\sfem}_n(y)=0$ for any $y\ge K$. Using
		$\Expect[p^\sfem_n(y)]=f_G(y)$ and, thus, $\EE[(p^\sfem_n(y) - f_G(y))^2] =
		\Var(p^\sfem_n(y))=\frac
		1{n^2}\sum_{i=1}^n\Var(\indc{Y_i=y})={f_G(y)(1-f_G(y))\over n}$ we continue the last display to get
		\begin{align}
			\frac{1}{c_2} \EE\qth{\dist(p^\sfem_n\|f_{G})\indc{Y_{\max}< K}}
			&\leq
			\sum_{y < K}{\EE\qth{(p^\sfem_n(y)-f_{G}(y))^2\indc{Y_{\max}< K}}\over f_{G}(y)}  
			+ \epsilon_K(G) (1-\epsilon_K(G))^n
			\nonumber\\
			&\leq
			\sum_{y < K}{\EE\qth{(p^\sfem_n(y)-f_{G}(y))^2}\over f_{G}(y)}  
			+ \epsilon_K(G) (1-\epsilon_K(G))^n
			\nonumber\\
			&\stepb{=}
			\frac{1}{n}\sum_{y< K} (1-f_{G}(y))
			+ \epsilon_K(G) (1-\epsilon_K(G))^n
			\leq {K\over n}+\epsilon_K(G). \label{eq:EB_m66}
		\end{align}
		Using the union bound and the fact $H^2\leq 2$ we have 
		$\EE\qth{H^2(f_{G},f_{{\hatGn}})\indc{Y_{\max}\geq K}}\leq 2\PP\qth{Y_{\max}\geq K}\leq 2n\epsilon_K(G).$
		Combining this with \prettyref{eq:EB_m6} and \prettyref{eq:EB_m66} yields
		\begin{align*}
			\EE\qth{H^2(f_G,f_{{\hatGn}})}
			&\leq \EE\qth{H^2(f_{G},f_{{\hatGn}})\indc{Y_{\max}< K}}
			+\EE\qth{H^2(f_{G},f_{{\hatGn}})\indc{Y_{\max}\geq K}}
			\nonumber\\
			&\leq \frac 4{c_1}\EE\qth{\dist(p^\sfem_n\|f_{G})\indc{Y_{\max}< K}}+2n\epsilon_K(G)
			\leq {4c_2\over c_1}{K\over n}+\pth{{4c_2\over c_1}+2n}\epsilon_K(G)\,,
		\end{align*}
		which completes the proof of~\eqref{eq:EB_general-H2-bound}. 
		
		To complete the proof of the theorem we need to choose the value of $K$ such
		that $\epsilon_K(G) \lesssim {1\over n^2}$. This is done slightly differently for
		each of the two different classes of priors:
		
		\begin{enumerate}[label=(\alph*)~]
			\item Pick $K=\min\sth{\ceil{5(he^2+2)\log n\over \log\log n}, he^2+5\log n}$. Then using \prettyref{lmm:bounded_prior_properties} with $a=4$ we get the desired result.
			\item 
			Let $G\in \subexpo(s)$.
			Choose $K=\max\sth{1,{2\over \log\pth{1+\frac 1{2s}}}}\log n$. Then \prettyref{lmm:subexpo_properties} in \prettyref{app:EB_subexpo-properties} implies that $\epsilon_K(G)\leq \frac 3{n^2}$. Plugging this in \eqref{eq:EB_general-H2-bound} completes the proof.
		\end{enumerate}
	\end{proof}
	
	\section{Proof of regret upper bound}
	\label{sec:EB_regret-upper-bound}
	\subsection{General regret upper bound via density estimation}
	The proof of \prettyref{thm:EB_main} relies on relating the regret in EB estimation to estimating the mixture density in the Hellinger distance. This idea has been previously noted in \cite[Theorem 3]{JZ09} for the Gaussian location models using Fourier analysis and an ingenious induction argument. Here the analysis turns out to be much simpler 
	thanks in part to the discreteness of the Poisson model and the light tail of the prior, leading to the following deterministic result which is crucial for proving the regret optimality of minimum-distance EB estimators.

	\begin{lemma}\label{lmm:EB_general-regret-bound}
		Let $G$ be a distribution such that $\EE_{G}[\theta^4]\leq M$ for some constant $M$. Then for any distribution $\hat {G}$ supported on $[0,\hat h]$, any $h>0$ with $G([0,h])>\frac 12$ and any $K\geq 1$,
		\begin{align*}
			\Regret(\hat G;G)
			\leq  
			&\sth{12(h^2+\hat h^2)
				+48(h+\hat h)K}(H^2(f_{{G}},f_{\hat{G}})+4G((h,\infty)))
			\nonumber\\
			&\quad +2(h+\hat h)^2\epsilon_K(G) +{2(1+2\sqrt 2)\sqrt{(M+\hat h^4)G((h,\infty))}}
		\end{align*}
		where 
		$\Regret(\hat G;G)$ and 
		$\epsilon_K(G)$ were defined in \eqref{eq:ep-K(G)} and \prettyref{eq:regret} respectively.
	\end{lemma}
	Note that   $\hat G$ in the above statement denotes, with a slight abuse of notation, an arbitrary deterministic prior and to compute the regret for the random estimator $\hat G$ we will need to consider the expectation. This will be clarified in the proof of \prettyref{thm:EB_main} later on. We provide a sketch of the proof here (see \prettyref{app:EB_general-regret-bound} for the full proof.) It is relatively easy to bound the regret if the corresponding Bayes estimator is also bounded, which is the case if the prior $G$ is compactly supported. Otherwise, one can consider its restriction $G_h$ on $[0,h]$ defined by $G_h(\cdot)={G(\cdot \cap[0,h])\over G([0,h])}$. The truncation error can be controlled  using properties of the mmse as follows:
	\begin{align}\label{eq:one-dim-subexpo-reduction}
		\Regret(\hat G;G)
		\leq  
		\Regret(\hat G;G_h)+{(1+2\sqrt 2)\sqrt{(M+\hat h^4)G((h,\infty))}\over G([0,h])}.
	\end{align}
	Then we use the structure of the Bayes estimator \prettyref{eq:EB_bayes_est} in the Poisson model to relate $\Regret(\hat G;G_h)$ to the squared Hellinger distance between $f_{G_h}$ and $f_{\hat G}$
	\begin{align}
		\Regret(\hat G;G_h)
		\leq \sth{6(h^2+\hat h^2)
			+24(h+\hat h)K}H^2(f_{{G}_h},f_{\hat{G}})
		+(h+\hat h)^2\epsilon_K(G_h),
		\label{eq:EB_m5}
	\end{align}
	for any $K\geq 0$. We then show that $\epsilon_K(G_h)$ and $H^2(f_{{G}_h},f_{\hat{G}})$ satisfies
	$$
	\epsilon_K(G_h)\leq 2\epsilon_K(G),\quad
	H^2(f_{\hat{G}},f_{{G}_h})
	\leq 2\sth{H^2(f_{G},f_{\hatGn})+4G((h,\infty))}.
	$$
	Replacing these bounds in \eqref{eq:EB_m5} we get the desired result.

	\subsection{Proof of \prettyref{thm:EB_main}}
	For rest of the section, let $C_1,C_2,\dots$ denote constants independent of $h,s,d$ as required. For Part (a), 
	recall that ${\hatGn}=\argmin_{Q\in \calP(\reals^+)} \dist\pth{p^\sfem_n\| f_{Q}}$ is the unconstrained minimum distance estimator.	
	To apply \prettyref{lmm:EB_general-regret-bound}, set $$\hat h=Y_{\max},\quad K=\min\sth{\ceil{5(he^2+2)\log n\over \log\log n}, he^2+5\log n},\quad M= h^4.$$ 
	Then, in view of \prettyref{lmm:bounded_prior_properties} we get that 
	\begin{align}\label{eq:r4}
	\PP\qth{\hat h> K}
	\leq n\cdot \PP\qth{Y_1> K}
	\leq {2\over n^4}.
	\end{align}
	For any $G\in\calP([0,h])$ we have from the proof of \prettyref{thm:EB_density-H2}(a)
	\begin{align*}
		\epsilon_K(G)\leq \frac {C_1}{n^3},\quad G((h,\infty))=0,\quad 
		\EE\qth{H^2(f_{G},f_{\hatGn})}
		\leq C_1\frac {K}n.
	\end{align*}
	Then \prettyref{lmm:EB_general-regret-bound} yields the required bound
	\begin{align*}
		\Regret({\hatGn};G)
		&\leq 
		C_3\sth{h^2+K^2+(h+K)K}\EE\qth{H^2\pth{f_{G},f_{\hatGn}}}+
		C_4\EE\qth{\pth{h^2+\hat h^2+(h+\hat h)K}\indc{\hat h>K}}
		\nonumber\\
		&+{\frac {C_4 \EE\qth{(h+\hat h)^2}}{n^3}}
		\nonumber\\
		&\stepa{\leq} C_5 \pth{\frac {K^3}n+\sqrt{\EE\qth{h^4+\hat h^4+K^4}\cdot \PP[\hat h>K]}}
		\stepb{\leq} {C_6K^3\over n},
	\end{align*}
	where (a) followed from the Cauchy-Schwarz inequality, and (b) from \eqref{eq:r4} and \prettyref{lmm:bounded_prior_properties}. 
	
	To achieve the result involving the constrained optimizer, we use $\hat h=h$. Following the proof of \prettyref{thm:EB_density-H2}, note that we only required the optimality property of the estimator of $G$ in \eqref{eq:EB_m6}, and the above equation holds true for the constrained estimator $\tilde G$ as well. Hence we can show 
	$$
	\sup_{G\in \calP[0,h]}\EE\qth{H^2(f_{G},f_{\tilde G})}\leq {c_1\over n}\cdot \min\sth{\max\{1,h\}{\log n\over \log\log n},h+\log n}
	$$ 
	as well. In view of \prettyref{lmm:EB_general-regret-bound} we get
	$
		\Regret({\hatGn};G)
		\leq 
		C_6\sth{h^2+hK}\EE\qth{H^2\pth{f_{G},f_{\tilde G}}}+{\frac {C_7 h^2}{n^3}}
		{\leq} C_8{hK^2\over n}.
	$



	For Part (b), we choose
	\begin{align}\label{eq:EB_m2}
		h={4s\log n},
		\quad K=\max\sth{1,{2\over \log\pth{1+\frac 1{2s}}}}\log n,
		\quad M=12s^4.
	\end{align}
	Since $G$ is $s$-subexponential, we have (see \prettyref{lmm:subexpo_properties} in \prettyref{app:EB_subexpo-properties} for details)
	\begin{align}\label{eq:EB_m3}
		\EE_G[\theta^4]\leq M,
		\quad G((h,\infty))\leq {2\over n^4},
		\quad \epsilon_K(G)\leq {3\over n^2},
		\quad \epsilon_{2K}(G)\leq {3\over n^4}.
	\end{align}
	In view of \prettyref{lmm:EB_support-g-hat} in \prettyref{app:EB_existence-uniqueness} we get that ${\hatGn}$ is supported on $[0,\hat h]$ where $\hat h=Y_{\max}$. Then \prettyref{lmm:EB_general-regret-bound} and $(\EE_G[Y_{\max}^2])^2\leq \EE_G[Y_{\max}^4]\leq \max\{1,s^4\}(\log n)^4$ (see \prettyref{app:EB_subexpo-properties} for a proof) implies
	\begin{align}
		\Regret({\hatGn};G)
		\leq 
		\EE\left[\sth{6(h^2+Y_{\max}^2)+24K(h+Y_{\max})}H^2(f_{G},f_{{\hatGn}})\right]+ {\frac {C_5}n}.
		\label{eq:EB_m8}
	\end{align}
	Next we bound the expectation in the last display. Using the fact that $H^2\leq 2$, we get 
	\begin{align}
		&\EE\left[\sth{(h^2+Y_{\max}^2)+4K(h+Y_{\max})}H^2(f_{{G}},f_{{\hatGn}}) \right]
		\nonumber\\
		&\leq  (h^2+4Kh+12K^2)\EE\qth{H^2\pth{f_{{G}},f_{\hatGn}}}
		+2\EE\qth{\sth{(h^2+Y_{\max}^2)+4K(h+Y_{\max})}\indc{Y_{\max}\geq 2K}}
		\label{eq:EB_m0}
	\end{align}
	Using \prettyref{thm:EB_density-H2} we get that the first part on the right of the above inequality is bounded as
	$$
		(h^2+4Kh+12K^2)\EE\qth{H^2\pth{f_{{G}},f_{\hatGn}}}
		\leq c_0\max\sth{1,s^3}{(\log n)^3\over n}
	$$
	for some absolute constant $c_0>0$. For the second term in \eqref{eq:EB_m0} we use Cauchy-Schwarz inequality and union bound to get
	\begin{align*}
		&\EE\qth{\sth{6(h^2+Y_{\max}^2)+24K(h+Y_{\max})}\indc{Y_{\max}\geq 2K}}
		\nonumber\\
		&{\leq} \sqrt{\EE\qth{\sth{6(h^2+Y_{\max}^2)+24K(h+Y_{\max})}^2}
			\PP_G\qth{{Y_{\max}\geq 2K}}}
		\nonumber\\
		&\stepa{\leq} \sqrt{\EE\qth{\sth{6(h^2+Y_{\max}^2)+24K(h+Y_{\max})}^2}
			n\epsilon_{2K}(G)}
		\nonumber\\
		&\stepb{\leq} {6\over n^{3/2}}\sqrt{\EE\qth{\sth{4(h^4+Y_{\max}^4)+16K^2(h^2+Y_{\max}^2)}}}
		\leq {1\over n}.
	\end{align*}
	where (a) followed from \eqref{eq:EB_m3} and (b) followed for large enough $n$. Plugging the bounds back in \eqref{eq:EB_m0} and in view of \prettyref{eq:EB_m8}, we complete the proof.

	\section{Numerical experiments}
	\label{sec:numerical-expt}
	In this section, we analyze the performances of the empirical Bayes estimators based on the minimum-$H^2$, the minimum-$\chi^2$, and the minimum-KL divergence estimator (i.e., the NPMLE). We compare them against the Robbins estimator and also draw comparisons among their individual performances. Unlike the Robbins estimator, the minimum-distance based estimators do not admit a closed form solution. 
	Our algorithm to compute the solution is closely related to the vertex direction method (VDM) algorithms for finding NPMLE \cite{L83general,L95}, specialized for the Poisson family and modified to work with the generalized distance we considered. In the case of the NPMLE, the convergence of the VDM method to the unique optimizer is well-known \cite{fedorov2013theory,wynn1970sequential}, and the algorithms for finding the other minimum $\dist$-distance estimators are expected to show similar convergence guarantees as well. Additionally, thanks to the Poisson density, the first-order optimality condition takes on a polynomial form, allowing us to use existing root-finding algorithms for polynomials to update the support points of the solution. 
	See \cite{S76} for a similar VDM-type algorithm for Poisson mixtures and \cite{KM14,koenker2017rebayes} for discretization-based algorithms.
	
	\subsection{First-order optimality condition and algorithm}
	\label{sec:first-order-optimality}
	In the numerical experiments we focus on the unconstrained minimum-distance estimator ${\hatGn}=\argmin_Q \dist(p^\sfem_n\|f_Q)$, which is a discrete distribution (\prettyref{thm:existence-uniqueness}).
	For any $\theta\in\reals_+$ let $\delta_\theta$ denote the Dirac measure at $\theta$. Suppose that the support of $p^\sfem_n$ be $\{y_1,\dots,y_m\}$. The optimality of $\hatGn$ implies that for all $\theta,\epsilon\in[0,1]$ we have $
		\dist(p^\sfem_n\|f_{{\hatGn}})\leq 	\dist(p^\sfem_n\|f_{(1-\epsilon){\hatGn}+\epsilon\delta_\theta}),$
	leading to the first-order optimality condition $\left.{d\over d\epsilon}\dist(p^\sfem_n\|f_{(1-\epsilon){\hatGn}+\epsilon\delta_\theta})\right|_{\epsilon=0}\geq 0$, namely
	\begin{align}\label{eq:EB_first-order-optimality}
		D_{\hatGn}(\theta) \triangleq \sum_{i=1}^m\left.\frac d{df}\ell(p^\sfem_n(y_i),f)\right|_{f=f_{\hatGn}(y_i)}
		(f_\theta(y_i)-f_{\hatGn}(y_i))\geq 0.
	\end{align}
	Averaging the left hand side over $\theta\sim\hatGn$, we get $\int D_{\hatGn}(\theta)d\hatGn(\theta)=0$. This implies that each $\theta$ in the support of $\hatGn$ satisfies $D_{\hatGn}(\theta)=0$. Taking derivative on both sides of the equation $D_{\hatGn}(\theta)=0$ with respect to $\theta$ we get that the atoms of ${\hatGn}$ satisfies the following polynomial equation in $\theta$
	\begin{align*}
		\sum_{i=1}^m w_i({\hatGn})\pth{y_i\theta^{y_i-1}-\theta^{y_i}}=0,\quad 
		w_i({\hatGn})=\sth{\left.\frac d{df}\ell(p^\sfem_n(y_i),f)\right|_{f=f_{{\hatGn}}(y_i)}}/ y_i!.
	\end{align*}
	Iterating the above conditions leads to following algorithm for computing $\hat G$.

	\begin{algorithm}[H]
		\caption{Computing the minimum $\dist$-distance estimators}\label{algo:proj-STG}
		{\bf Input}:  Data points $Y_1,\dots,Y_n$. Target distribution $G_{\boldsymbol\theta,\boldsymbol\mu}=\sum_{j}\mu_j\delta_{\theta_j}$. Divergence $\dist$ with $t-\ell$ decomposition $\dist(q_1\|q_2)=t(q_1)+\sum_{y\geq 0}\ell(q_1(y),q_2(y))$. 
		Initialization of $(\boldsymbol{\theta},\boldsymbol{\mu})$. Tolerance $\eta_1,\eta_2$ and number of iterations $N$. \\
		{\bf Steps}: 
		\begin{algorithmic}[1]
			\State Calculate empirical distribution $p_n^\sfem$. Obtain the set of distinct sample entries $\sth{y_1,\dots,y_m}$.
			\For{$N$ iterations }
			\State $newroots=\sth{\theta:\theta\geq 0, \sum_{i=1}^m w_i({ G}_{\boldsymbol\theta,\boldsymbol\mu})\pth{y_i\theta^{y_i-1}-\theta^{y_i}}=0}$.
			\State Combine $\boldsymbol\theta$ and $newroots$ and denote the new vector as $\boldsymbol\theta'$.
			\State Merge entries of $\boldsymbol\theta'$ that are within $\eta_1$ distance of each other.
			\State Find $\argmin_{\tilde{\boldsymbol\mu}}\sum_{i=1}^m\ell(p^\sfem_n(y_i),f_{G_{\boldsymbol\theta',\tilde{\boldsymbol\mu}}}(y_i)),$ via gradient descent with initialization $\tilde{\boldsymbol\mu}=\boldsymbol \mu$.
			\State Remove entries of $\boldsymbol\theta'$ and $\boldsymbol\mu'$ at locations of $\boldsymbol\mu'$ that are less than $\eta_2$ and re-normalize $\boldsymbol\mu'$.
			\State $(\boldsymbol\theta,\boldsymbol\mu)\leftarrow (\boldsymbol\theta',\boldsymbol\mu')$.
			\EndFor
		\end{algorithmic}
		{\bf Output}: $(\boldsymbol{\theta},\boldsymbol{\mu})$.
	\end{algorithm}
	We apply this algorithm for finding the minimum-distance estimators in the following examples. In all our experiments we used $\eta_1=0.01,\eta_2=0.001$. We set the maximum number of iterations $N$ to be 15 as the outputs of the algorithm in all our simulations were observed to converge by then.
	We choose the initialization for $\boldsymbol\theta$ as the uniform grid of size 1000 over the interval $[0,Y_{\max}]$, with a uniform initial probability assignment $\boldsymbol\mu$.

	\subsection{Real-data analysis: Prediction of hockey goals}
	\label{sec:EB_real_data}
\begin{figure}[t]
		\centering
		\begin{minipage}{0.5\textwidth}
			\begin{center}
				\includegraphics[height=6cm]{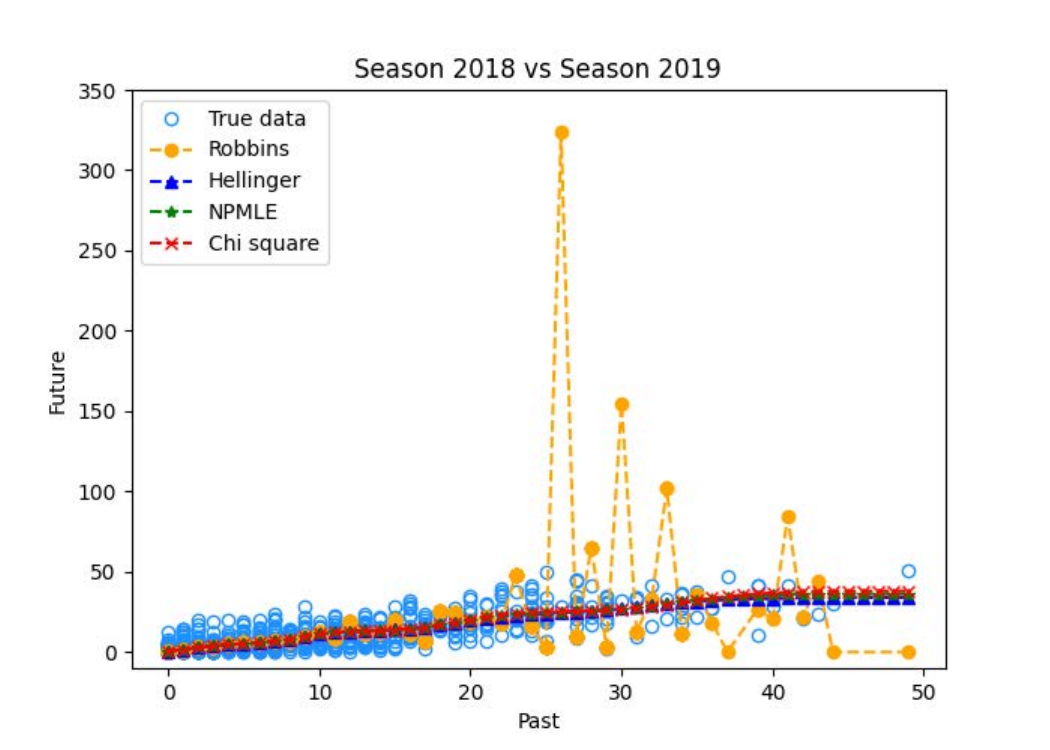}
			\end{center}
		\end{minipage}\hfill
		\begin{minipage}{0.5\textwidth}
			\begin{center}
				\includegraphics[height=6cm]{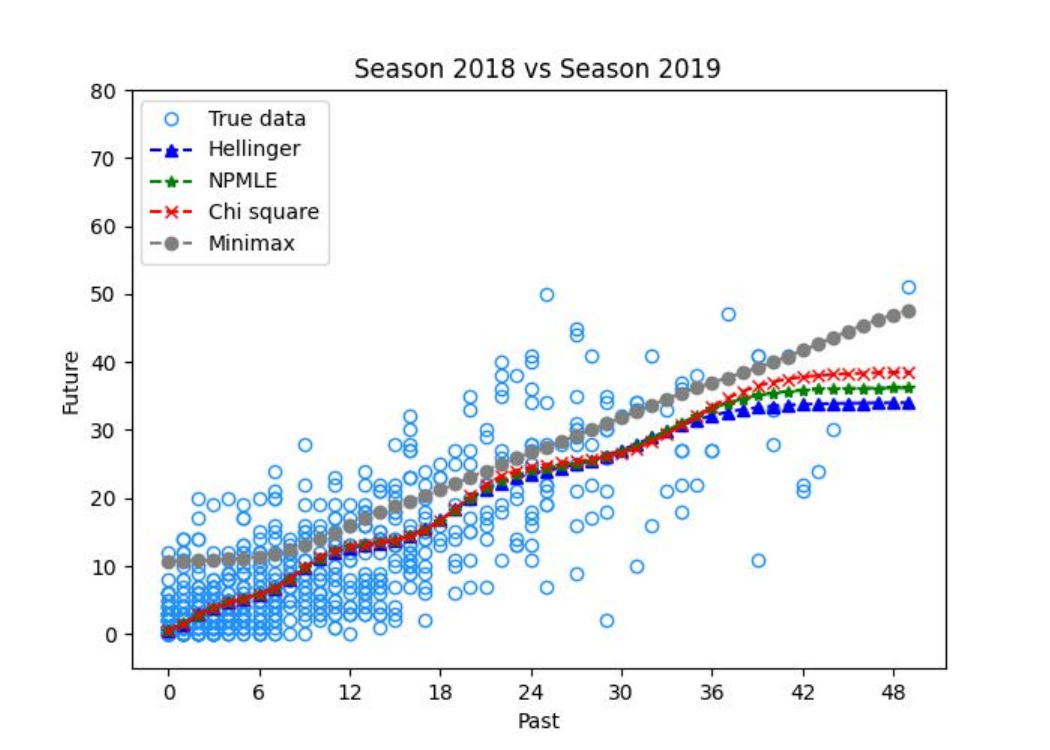}
			\end{center}		
		\end{minipage}
		\caption{Prediction of hockey goals with empirical Bayes, comparing Robbins and minimum-distance estimators. On the right panel, the Robbins estimator is replaced by the minimax estimator.}
		\label{fig:real}
	\end{figure}

	Here we extend our study related to the National Hockey League as mentioned in \prettyref{sec:EB_introduction}. 
	{  In \prettyref{fig:real} we first plot the actual data, where for each true data point, its $x$-axis represents the number of goals scored by a particular player in the 2017-18 season (denoted as ``Past") and the $y$-axis of the data represents the goals by the same player in the 2018-19 season (denoted as ``Future"). Then, for each possible value of the goal scored in the 2017-18 season, we plot the EB estimators based on the Robbins method, the minimum $H^2$, the minimum-$\chi^2$ distance estimator, and the NPMLE, and the {\it gold standard} minimax estimator. }
	
	The left panel shows that a large number of individuals exist for whom the Robbins estimator produces unstable predictions that are significantly worse than those of all minimum-distance methods. This difference is significant for the values of scored goals, which have lower sample representations. 
	Thus, on the right panel, we omit the Robbins estimator and provide a more detailed comparison of
	the three minimum-distance estimators, which shows that their behavior is mostly comparable
	except near the tail end of the data points. We also present the comparison with the minimax estimator given by the conditional estimator of the least favorable prior in \prettyref{fig:least_fav_prior}, and our plot shows that the EB estimators align more closely to the actual data. 
	\begin{figure}[t]
		\centering
		\begin{minipage}{0.5\textwidth}
			\begin{center}
				\includegraphics[ height=6cm]{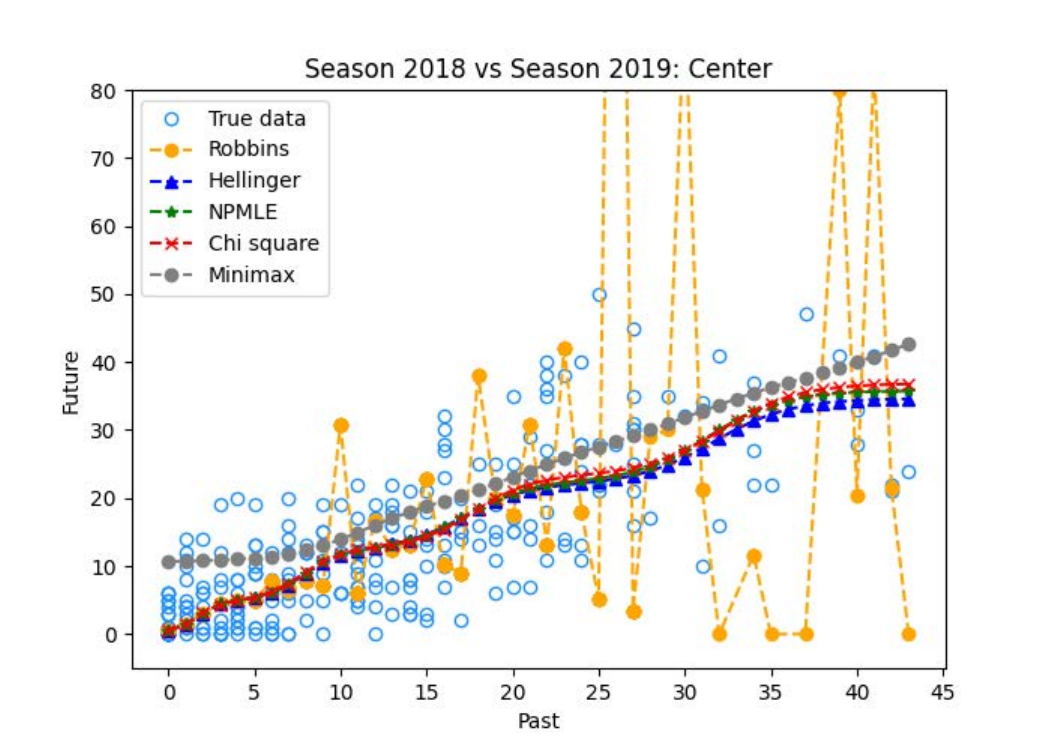}
			\end{center}
		\end{minipage}\hfill
		\begin{minipage}{0.5\textwidth}
			\begin{center}
				\includegraphics[ height=6cm]{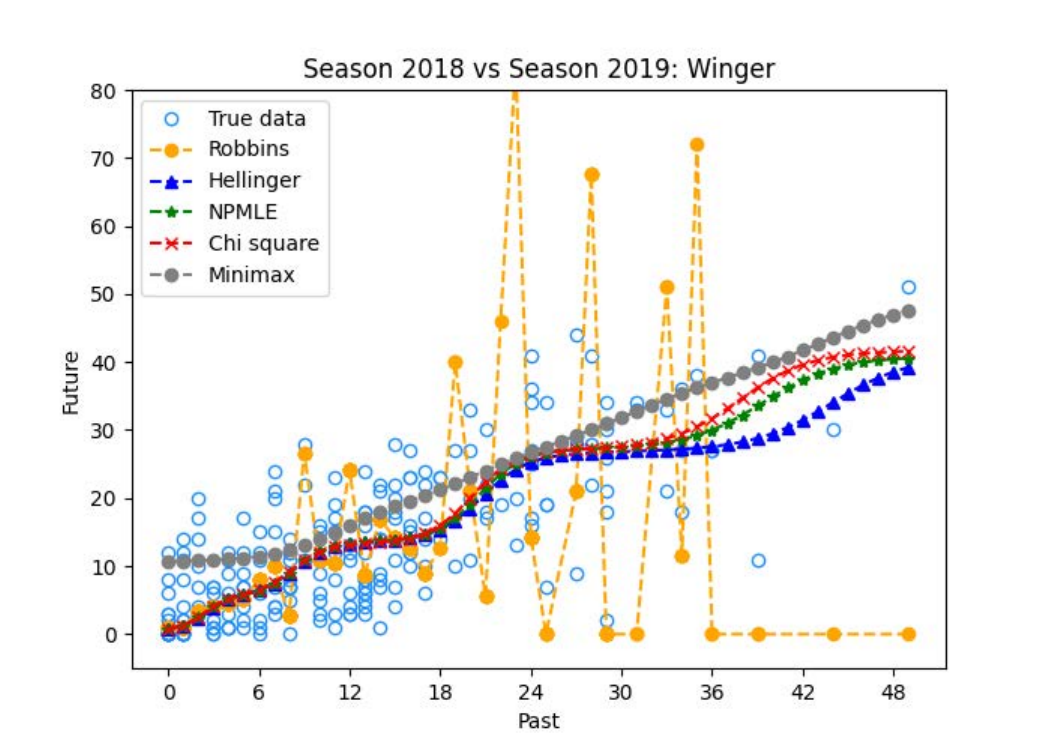}
			\end{center}		
		\end{minipage}
		\begin{minipage}{0.5\textwidth}
			\begin{center}
				\includegraphics[ height=6cm]{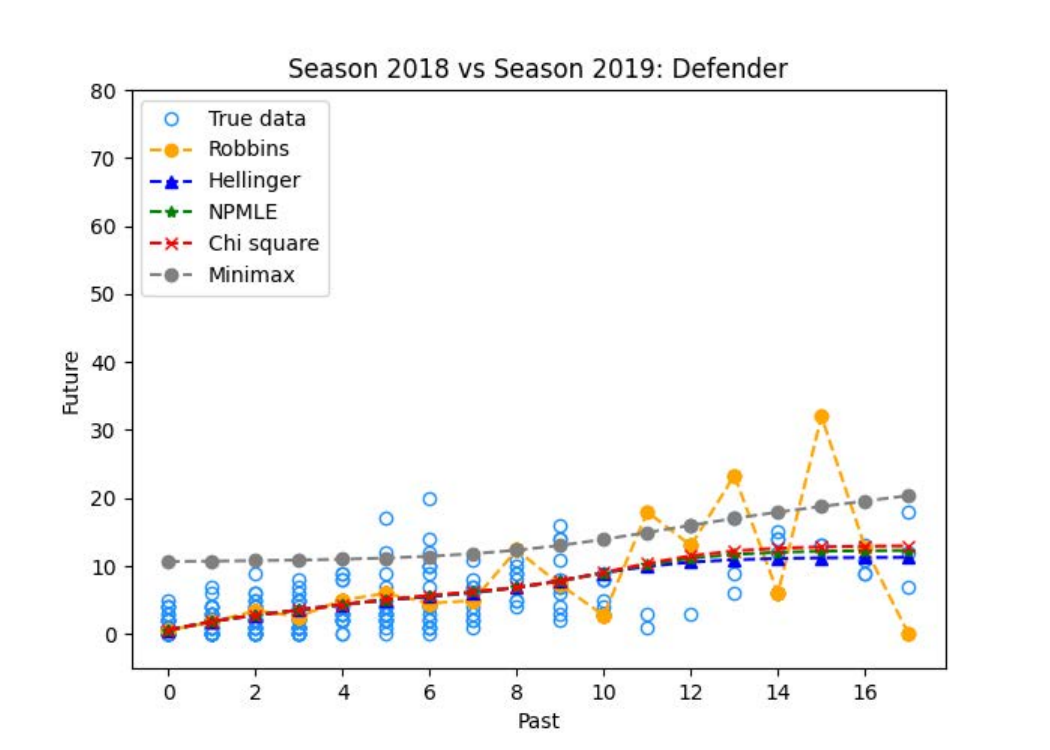}
			\end{center}		
		\end{minipage}
		\caption{Prediction of hockey goals at different playing positions.}
		\label{fig:real_position}
	\end{figure}
	
	Interestingly, all three estimators seem to do shrinkage towards several fixed values. There could be several explanations for this multi-modality. One is that different clusters correspond to different player positions (defense, winger, center). The other is that clusters correspond to the line of the player (different lines get different amounts of ice time). To test this hypothesis, we also redid on Fig. 3 the estimation for each position separately. Since the multi-modality is retained, we conclude that the second option is more likely to be the real explanation. In addition, we also compared the four goal-prediction methods based on different EB estimators and the minimax estimator based on the least favorable prior on $[0,50]$ across the possible playing positions: defender, center, and winger. Similar to before, we used the Poisson model and tried to predict the goal scoring for the year 2019 using the goal scoring data from the year 2018 for players in each playing position separately. As expected, the minimum distance methodology provides more stable and accurate estimates than the estimates based on the Robbins method and minimax strategy. The plots showing the closeness of the predictions to the actual number of goals for the different EB methods are provided in \prettyref{fig:real_position}.
	
	\subsection{Application of EB methods for filtering}
	In this section, we demonstrate the application of the EB methodology to data cleaning. We propose to show that given a standard statistical methodology, incorporating an EB-based filter on the data before feeding it to the algorithm, can significantly improve the existing performance guarantees. For our analysis we use simulated data based on multivariate linear models.
	
	For this simulation study, we assume that the observed data $y_1,\dots,y_n$ are independently generated via a linear model
	$$
	y_i \overset{\text{ind.}}{=} {\boldsymbol\theta}_i{\beta},\quad {\boldsymbol\theta}_i {=\sth{\theta_{ij}}_{j=1}^d}\in \reals_+^d, \beta\in \reals^d, \quad i=1,\dots,n.
	$$ 
	In addition, we also assume that the observer does not directly see the data generating $\theta_i$-s. Instead, we can only observe a Poissonized version $X_i$ of $\theta_i$, given as
	$$
	{\boldsymbol X}_i=\{X_{ij}\}_{j=1}^d,\quad X_{ij}\inddistr \Poi(\theta_{ij}),\quad i=1,\dots,n,\quad j=1,\dots,d.
	$$  
	In other words, each coordinate of $\boldsymbol X_i$ is generated independently according to a Poisson channel, with mean being the corresponding coordinate of $\boldsymbol \theta_i$. Then we pose the following question: 
	\begin{align*}
		&\text{\it   Upon observing $\sth{(y_i,{\boldsymbol X}_i)}_{i=1}^n$, can we achieve better error guarantees if we apply the EB filters}\\ 
		&\text{\it   to the covariates before running the ordinary least squares (OLS) methodology?}
	\end{align*}
	
	\begin{figure}[t]
		\centering
		\includegraphics[height=8cm]{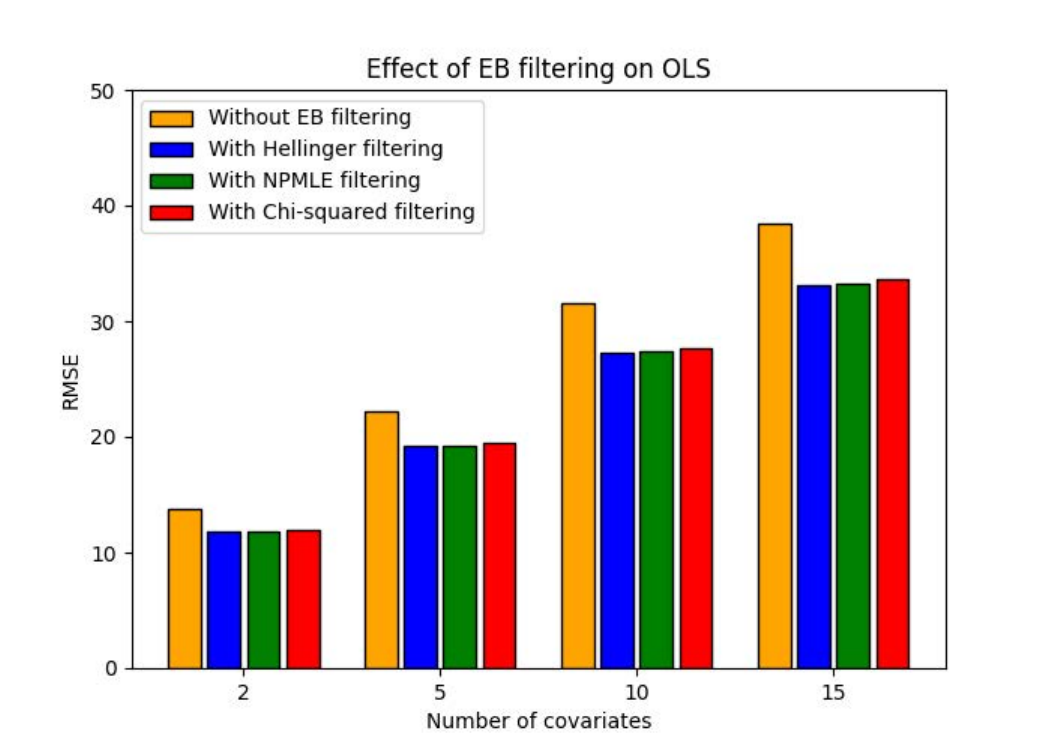}
		\caption{Improving results for OLS using EB filtering}
		\label{fig:regression}
	\end{figure}
	
	To answer the problem, we will show that EB-based one-dimensional data filters, applied separately on each of the covariates before running the OLS, can significantly improve the prediction of $y$ if the ${\boldsymbol\theta}_i$-s values are coming from a multivariate mixture model. Our process for generating the data $\sth{(y_i,{\boldsymbol X}_i)}_{i=1}^n$ is outlined as follows. To generate each coordinate of ${\boldsymbol\theta}_i$, we randomly generated entries from a uniform mixture of Gaussian distributions, with means $[2,8,16,32]$ and standard deviation 1, and then took the absolute values. Each coordinate of the regression coefficient $\beta$ was selected uniformly from $[-5,5]$. The number of covariates (i.e., $d$) are varied in the set $[2,5,10,15]$. A total of 1200 samples $(y_1,...,y_{1200})$ were generated. To assess the goodness of the fit, we compute $\hat y$ and the corresponding root mean squared error (RMSE). 
	To examine the effect of the one-dimensional EB filters on the covariates, we used NPMLE, Hellinger-based, and Chi-squared distance based one dimensional EB filters. The RMSE prediction errors were compared with and without the filtering. In all the simulations, the EB filtering improved the result. The plot of errors is presented in \prettyref{fig:regression}, and the errors are reported (along with standard deviations inside parenthesis) in \prettyref{tab:regression}.

	\begin{table}[h]
		\centering
		\caption{Performance of EB filtering}
		\label{tab:regression}
		\begin{tabular}{@{\extracolsep\fill}lllll@{\extracolsep\fill}}
			\hline
			Covariates & Without filtering & $H^2$ filter & NPMLE filter & $\chi^2$ filter\\
			\hline
			2 & 13.75 (0.149) & 11.857 (0.131) & 11.867 (0.131) & 11.993 (0.132)\\
			\hline
			5 & 22.25 (0.145) & 19.262 (0.128) & 19.277 (0.128) & 19.497 (0.129)\\
			\hline
			10 & 31.621 (0.143) & 27.327 (0.126) & 27. 349 (0.126) & 27.670 (0.128)\\
			\hline
			15 & 38.412 (0.145) & 33.163 (0.129) & 33.195 (0.129) & 33.578 (0.130)
			\\
			\hline
		\end{tabular}
	\end{table}

	\subsection{More simulation studies}
	\label{sec:simulation-unbounded}
	In this subsection, we test more priors in addition to the uniform prior in \prettyref{fig:Rob-vs-NPMLE}, including discrete priors and priors with unbounded support. In \prettyref{sec:EB_real_data} we see that the three minimum-distance estimators performed similarly. However, the question arises whether the best choice among the minimum-distance EB methods can be argued when some information about the prior is available. With the specific goal of differentiating the three minimum-distance estimators among themselves, we carry out simulation studies at the end of this section using different priors.
	
	\begin{figure}[t]
		\centering
		\begin{minipage}{0.5\textwidth}
			\begin{center}
				\includegraphics[ height=6cm]{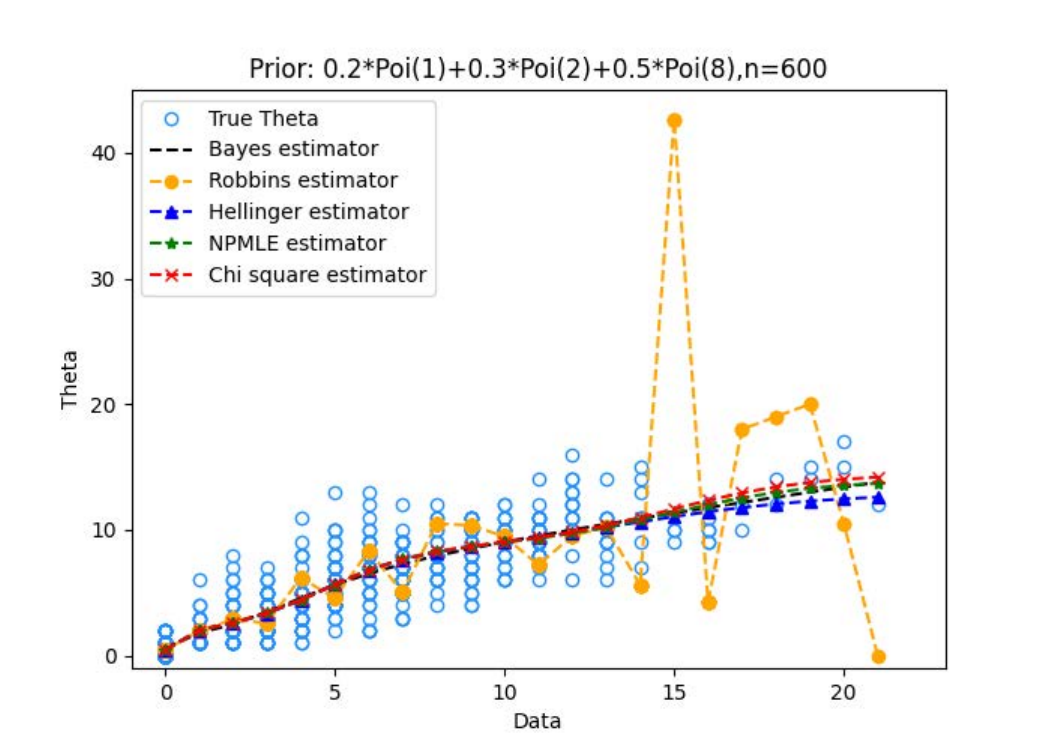}
			\end{center}
		\end{minipage}\hfill
		\begin{minipage}{0.5\textwidth}
			\begin{center}
				\includegraphics[ height=6cm]{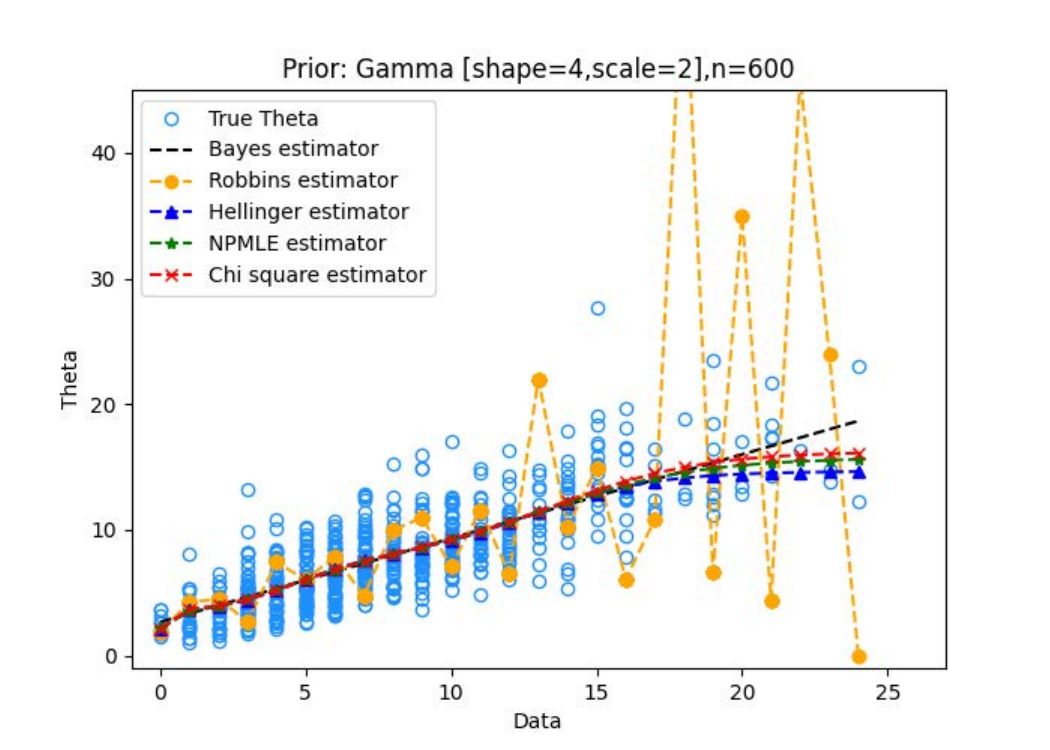}
			\end{center}
		\end{minipage}
		\caption{Robbins vs. minimum-distance estimators: Unbounded priors}
		\label{fig:unbounded-priors}
	\end{figure}
	
	For comparing the EB methods in the discrete setup we choose the prior $G$ to be $0.2 \Poi(1)+0.3\Poi(2)+0.5\Poi(8)$ and for the continuous unbounded setup we choose the prior $G$ to be the Gamma distribution with scale parameter 2 and shape parameter 4, i.e., with prior density $f(x)=\frac 1{96}x^3e^{-\frac x2}$. In both of the cases we simulate $\sth{\theta_i}_{i=1}^{600}$ independently from the prior distribution and correspondingly generate data $Y_i\sim\Poi(\theta_i)$. For each of the priors we calculate the Bayes estimator numerically (denoted by the black dashed line in the plots). Then, from the generated datasets, we compute the Robbins estimator, the NPMLE-based EB estimator, the $ H^2$-distance-based EB estimator, and the $\chi^2$-distance-based EB estimator. All the estimators are then plotted against $\theta$ and the data (\prettyref{fig:unbounded-priors}). As expected, the Robbins estimator shows high deviation from the true $\theta$ values in many instances whereas the minimum-distance based estimators are much more stable.
	
\begin{figure}[t]
		\centering
		\begin{minipage}{0.5\textwidth}
			\begin{center}
				\includegraphics[ height=6cm]{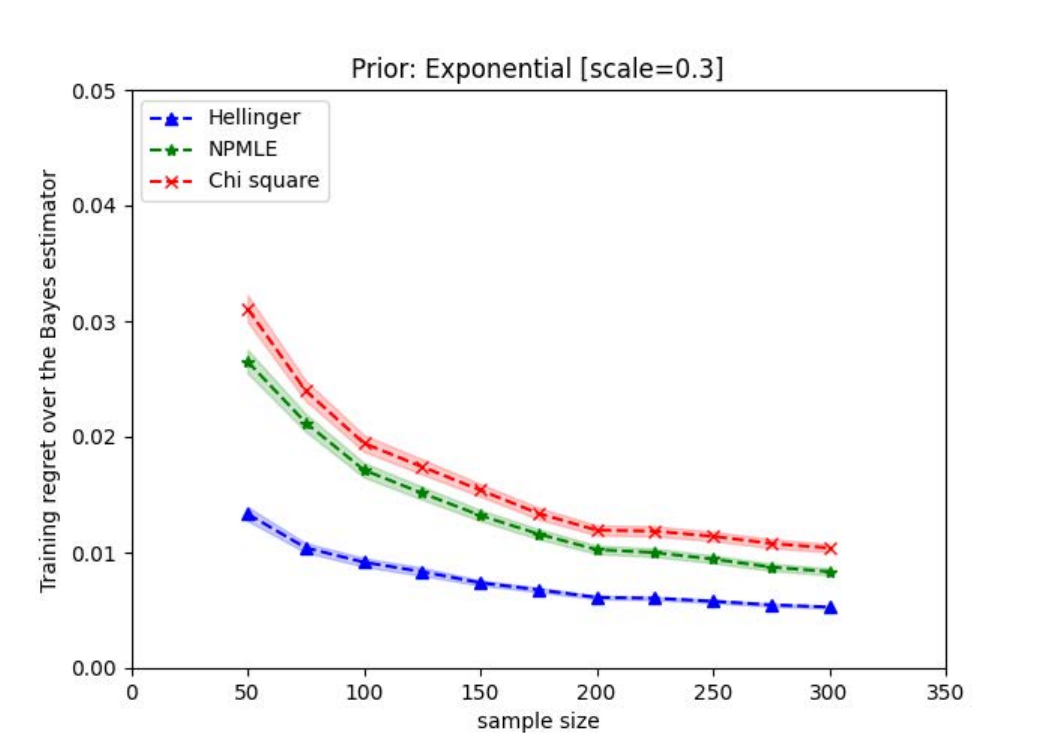}
			\end{center}
		\end{minipage}\hfill
		\begin{minipage}{0.5\textwidth}
			\begin{center}
				\includegraphics[ height=6cm]{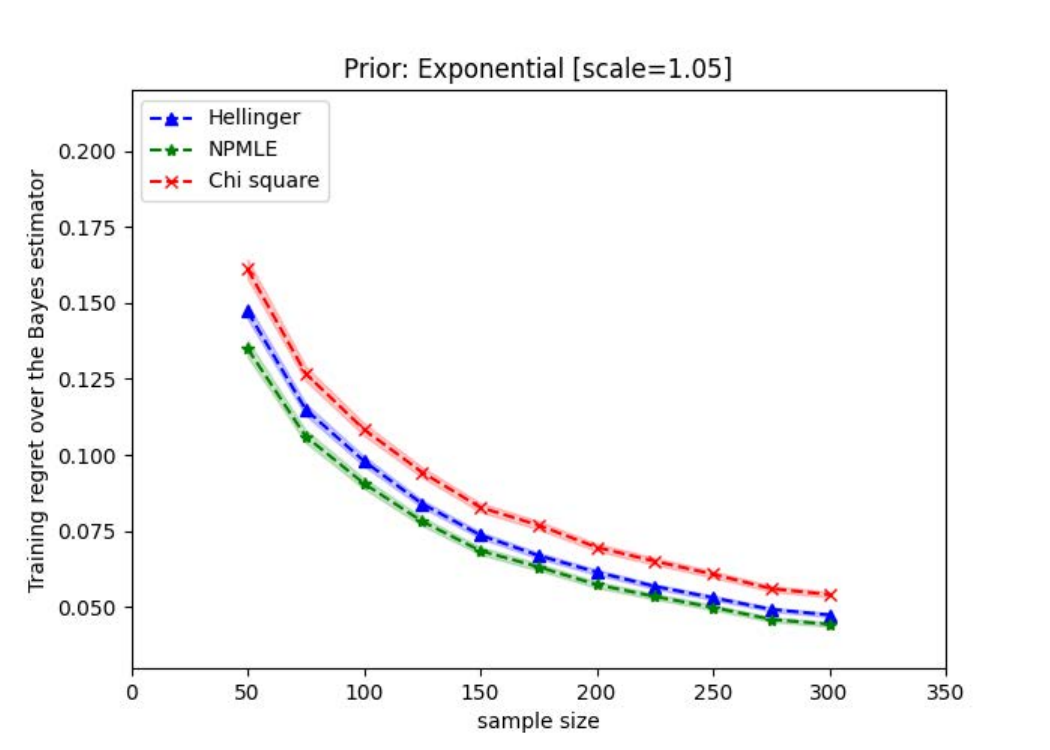}
			\end{center}
		\end{minipage}
		\begin{center}
			\includegraphics[ height=6cm]{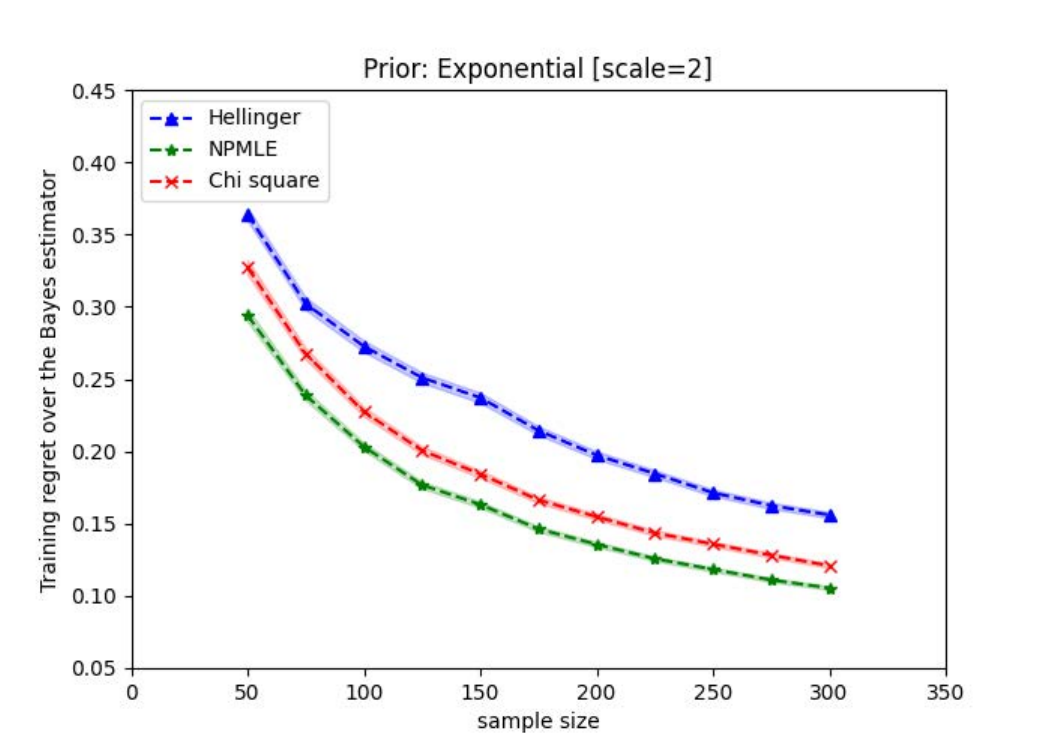}
		\end{center}
		\caption{Comparison of minimum-distance estimators}
		\label{fig:min-dist-comparison}
	\end{figure}

	To differentiate the different minimum-distance based EB methods we analyze the effect of the tail properties of the prior in the simulations below. Consider the exponential distribution parameterized by scale ($\alpha$) and with density $g_\alpha(x)=\frac 1\alpha e^{-x/\alpha}$. Note that the higher values of $\alpha$ generate distributions with heavier tails. We consider three values of $\alpha$: 0.3,1.05 and 2. For each $\alpha$ we estimate the training regret for sample sizes $n$ in the range $[50,300]$. Given sample $Y_1,\dots,Y_n$ from the mixture distribution with prior $G$ we define the training regret for any estimator $\hatGn$ of $G$ as $\EE_G[\frac 1n\sum_{i=1}^n(\hat\theta_G(Y_i)-\hat\theta_{\hatGn}(Y_i))^2]$. We compute the Bayes estimator $\hat \theta_G(y)$ numerically for each $y$. For every pair $(\alpha,n)$ we replicate the following experiment independently 10,000 times for each minimum-distance method: 
	
	\begin{itemize}
		\item ~Generate $\{\theta_i\}_{i=1}^n$ and $Y_i\sim\Poi(\theta_i)$, 
		\item ~Calculate ${\hatGn}$ using minimum-distance method,
		\item ~Calculate prediction error $\textsf{E}(Y^n)=\frac 1n\sum_{i=1}^n(\hat \theta_G(Y_i)-\hat\theta_{\hatGn}(Y_i))^2$.
	\end{itemize}
	Then we take the average of $\textsf{E}(Y^n)$ values from all the 10,000 replications to estimate the training error. For each $\alpha$ and each minimum distance method, at every $n$ we also estimate the $95\%$ confidence interval as $[\overline{\textsf{E}(Y^n)}\mp 0.0196 * \textsf{sd}(\textsf{E}(Y^n))]$ where $\overline{\textsf{E}(Y^n)}$ and $\textsf{sd}(\textsf{E}(Y^n))$ define respectively the sample mean and the sample standard deviation of the $\textsf{E}(Y^n)$ values over the 10,000 independent runs. Below we plot the training regrets and their 95\% confidence bands against the training sample sizes (\prettyref{fig:min-dist-comparison}). We observe that that minimum-$H^2$ based estimator outperforms the other estimators when the scale of the exponential distribution is small. As the tails of the prior distributions become heavier, the performance of the minimum-$H^2$ based estimator gets worse and the NPMLE based estimator comes out as a better choice.

	\section{An extension of the results in multiple dimensions}
	\label{sec:results-multidim}
	For a clarity of notations, in this section we use the bold fonts to denote vectors, e.g., $\btheta=\pth{\theta_1,\dots,\theta_d},\vtheta_{i}=(\theta_{i1},\dots,\theta_{id}),\vY=(Y_1,\dots,Y_d),\vY_{i}=(Y_{i1},\dots,Y_{id}), \vy=(y_1,\dots,y_d)$, etc. Let $G$ be a prior distribution on $\reals_+^d$. We will study an extension of our minimum distance estimators \eqref{eq:mindist} in the $d$-dimension coordinate-wise independent Poisson model: Consider the following data-generating process \cite{johnstone1986admissible,brown1985complete}
	\begin{align}
		\vtheta_i\iiddistr G, \quad Y_{ij} \inddistr \Poi(\theta_{ij}),\quad 
		i=1,\dots,n,\quad j=1,\dots,d.
	\end{align}
	Note that the marginal distribution of the multidimensional Poisson mixture is given by 
	$$f_G(\vy)=\int_{\vtheta} \prod_{i=1}^d e^{-\theta_i}{\theta_i^{y_i}\over y_i!} dG(\vtheta),\quad \vy\in\integers_+^d.$$
	To construct the minimum distance estimator we use the same minimization principle as in \eqref{eq:mindist}, where we will specify the prior class $\calG$ used for optimization when we state the results. Next we construct the EB estimator. Denote by 
	$\hat \vtheta_G$ the Bayes estimator, whose $j$-th coordinate $\hat \theta_{G,j}$ is
	$$
	\hat \theta_{G,j}(\vy)
	=\EE_G[\theta_{j}| \vy]
	={\int_{\vtheta} \theta_j \prod_{j=1}^d e^{-\theta_j}{\theta_j^{y_j}\over y_j!} dG(\vtheta)\over f_{G}(\vy)}=(y_j + 1)\frac{f_{G}(\vy + \ve_j)}{f_{G}(\vy)},\quad j=1,\dots,d,$$
	where $\ve_j$ denote the $j$-th coordinate vector. Suppose that $\hat G$ gives us an estimate of the prior distribution $G$
	and consider the corresponding empirical Bayes estimator 
	$\hat \vtheta_{\hat G} = (\hat \theta_{\hat G,1},\dots,\hat \theta_{\hat G,d}). $
	Similar to \prettyref{eq:regret}, let us define the regret of any plug-in estimator based on a prior estimate $\hat G$ as
	\begin{align}
		\Regret(\hat G,G)
		={\EE_G\qth{\|\hat\vtheta_{{\hatGn}}(\vY)-\hat\vtheta_{ G}(\vY)\|^2}} 
		= \EE_G\qth{\sum_{\vy\in \integers_+^d}
			\|\hat\vtheta_{{\hatGn}}(\vy)-\hat\vtheta_{ G}(\vy)\|^2 f_G(\vy)},
	\end{align}
	where $\vY\sim f_G$ is a test point independent from the training sample $\vY_1,\ldots,\vY_n\iiddistr f_G$ ($\integers_+^d$ denotes the set of all $d$-dimensional vectors with non-negative integer coordinates). We will prove regret bounds for the minimum distance estimator of the form \eqref{eq:mindist} where the $\dist$ function satisfies the following regularity assumption.
	
	\begin{assumption}
		\label{pt:EB_multi_t-ell} 
		There exist maps $\vt:\calP(\integers_+^d)\to \reals$ and $ \ell:\reals^2 \to\reals$ such that for any two distributions $q_1,q_2\in\calP(\integers_+^d)$
		$$\dist\pth{q_1\| q_2}= \vt(q_1)+\sum_{\vy\in \integers_+^d}\ell(q_1(\vy),q_2(\vy)),$$
		where $b \mapsto \ell(a,b)$ is strictly decreasing and strictly convex for $a> 0$ and $\ell(0,b)=0$ for $b\geq 0$. 
	\end{assumption}
	\begin{assumption}
		\label{pt:EB_multi_sandwich}
		There exist absolute constants $c_1,c_2>0$ such that for pmf-s $q_1,q_2$ on $\integers_+^d$
		\begin{align}\label{eq:EB_multi_Assupmtion_2}
			c_1H^2(q_1,q_2)\leq \dist(q_1\|q_2)\leq c_2\chi^2(q_1\|q_2).
		\end{align}
	\end{assumption}
	\prettyref{pt:EB_multi_t-ell} and \prettyref{pt:EB_multi_sandwich} are identical to \prettyref{pt:EB_t-ell} and \prettyref{pt:EB_sandwich} respectively. The loss functions such as the Kullback-Leibler divergence, squared Hellinger distance, Chi-squared divergence satisfy the above assumptions, similarly as in the one-dimensional case. We have the following results.
	
	\begin{theorem}\label{thm:density_multidim}
		Let $\dist$ satisfy \prettyref{pt:EB_multi_t-ell} and \prettyref{pt:EB_multi_sandwich}. Suppose that $\hat G$ is the unconstrained minimum distance estimator
		\begin{align}
			\label{eq:multidim-est-unconstrained}
			\hat G = \argmin_{Q\in\calP(\reals_+^d)} \dist(p_n^\sfem\|f_{Q}).
		\end{align}
		Then there exist constants $c_1,c_2$ such that the following holds
		\begin{enumerate}[label=(\roman*)~]
			\item $\sup_{G\in \calP([0,h]^d)}\EE\qth{H^2(f_{\hat G},f_G)}\leq \frac {\pth{c_1 K}^{d}}n$, where $K=\min\sth{\max\{1,h\}{\log n\over \log\log n},h+\log n}$ ; 
			
			\item If the data generating prior $G$ belongs to a class $\calG$ where all marginals of $G\in\calG$ belong to the $\subexpo(s)$ class of distributions
			for some $s>0$, then 
			$$\sup_{G\in \calG}\EE\qth{H^2(f_{\hat G},f_G)}\leq \frac {\pth{c_2 \max\{1,s\}}^{d}}n (\log (n))^{d}.$$
		\end{enumerate}
	\end{theorem}
	
	\begin{theorem}\label{thm:main_multidim}
		Suppose that the assumptions in \prettyref{thm:density_multidim} hold true. Then the following regret bounds hold for the unconstrained estimator \eqref{eq:multidim-est-unconstrained} whenever $n \ge d$ ($c_1, c_2> 0$ below are constants): 
		\begin{enumerate}[label=(\roman*)~]
			\item $\Regret(\hat G;\calP([0,h]^d))\leq \frac {d\pth{c_1 K}^{d+2}}n$, where $K=\min\sth{\max\{1,h\}{\log n\over \log\log n},h+\log n}$ ; 
			
			\item If the data generating prior $G$ belongs to a class $\calG$ where all marginals of $G\in\calG$ belong to the $\subexpo(s)$ class of distributions
			for some $s>0$, then 
			$$\Regret(\hat G;\calG)\leq \frac {d\pth{c_2 \max\{1,s\}}^{d+2}}n (\log (n))^{d+2}.$$
		\end{enumerate} 
		In addition, in the case when the data generating distribution $G$ is supported on $[0,h]^d$, then the constrained minimum distance estimator with access $h$ achieves improved risk guarantee
		\begin{align*}
			\tilde G=\argmin_{Q\in \calP[0, h]^d}\dist(p_n^\sfem\|f_{Q}),
			\quad
			\Regret({\tilde G};\calP([0,h]^d))
			\leq \frac {dc_1^{d+2} {\max\{1,h\}}}n {K}^{d+1}.
		\end{align*}
	\end{theorem}
	The proofs of the above results are provided in \prettyref{app:multidim-proofs} below. We conjecture these regret bounds in Theorem \ref{thm:main_multidim} are nearly optimal and factors like $(\log n)^d$ are necessary. A rigorous proof of matching lower bound for Theorem \ref{thm:main_multidim} will likely involve extending the regret lower bound based on Bessel kernels in \cite{PW21}
	to multiple dimensions; this is left for future work.
	
	\section*{Data availability}
	The real data set on hockey goals that we used is available at \href{https://www.hockey-reference.com/}{https://www.hockey-reference.com/}. Our code is available at \href{https://github.com/janasoham/codes_public/tree/main/mindist_poisson_eb}{https://github.com/janasoham/codes\_public/tree/main/mindist\_poisson\_eb}.
	
	\section*{Acknowledgment}
	Y.~Polyanskiy is
	supported in part by the MIT-IBM Watson AI Lab, and the NSF Grants CCF-1717842, CCF-2131115. Y.~Wu is supported in part by the NSF
	Grant CCF-1900507, NSF CAREER award CCF-1651588, and an Alfred Sloan fellowship.
		
	\bibliographystyle{apalike}
	\bibliography{refs_EB}
	
	\appendix

	{  \section{Least favorable prior and minimax estimator in the Poisson setup}
		
	\label{app:least-favorable}
	\begin{theorem}\label{thm:least-favorable}
		Let $h>0$ be finite and $\calP([0,h])$ denote the set of all probability distribution supported on $[0,h]$. Consider the minimax objective 
		$$
		\inf_{\hat \theta} \sup_{\theta \in [0,h]} \EE\qth{(\theta - \hat \theta(Y))^2},
		$$
		where $Y|\theta\sim \Poi(\theta)$. Then the least favorable prior is a discrete distribution. In addition, the minimax estimator for the objective 
		$$\inf_{\hat\theta_1,\dots,\hat \theta_n} \sup_{\theta_1,\dots,\theta_n\in [0,h]} \EE\qth{\sum_{i=1}^n(\theta_i - \hat \theta_i(Y^n))^2},
		\quad 
		Y^n=\{Y_i\}_{i=1}^n,Y_i|\theta_i\inddistr \Poi(\theta_i),i\in [n],
		$$ 
		is the conditional mean of the above least favorable distribution given $\{Y_1,\dots,Y_n\}$. 
	\end{theorem}
	\begin{proof}
	We will use results from \cite{dytso2018structure} with the notations
	$$
	\theta=X,\quad P_{Y|\theta}=\Poi(\theta).
	$$ 
	The discreteness of the least favorable prior follows from \cite[Proposition 6]{dytso2018structure}. Let $\hat G$ be the least favorable prior and $\hat \theta_{\hat G}(Y)$ be the conditional mean of the least favorable prior evaluated at $Y$. We show that $\sth{\hat\theta_{\hat G}(Y_i)}_{i=1}^n$ is the minimax estimator. Note that in view of \cite[Theorem 6]{dytso2018structure} we get
	$$
		\EE_{G}\qth{\pth{\theta-\hat \theta_{\hat G}(Y)}^2}
		\leq \EE_{\hat G}\qth{\pth{\theta-\hat \theta_{\hat G}(Y)}^2},
		\quad 
		G\in \calP([0,h]),Y|\theta\sim \Poi(\theta).
	$$
	This implies 
	$$
		\sup_{\theta_1,\dots,\theta_n} \sum_{i=1}^n
		\qth{\pth{\theta_i-\hat \theta_{\hat G}(Y_i)}^2}
		\leq 
		\EE_{\theta_i\sim \hat G,i\in[n]}\sum_{i=1}^n\qth{\pth{\theta_i-\hat \theta_{\hat G}(Y_i)}^2},
		\quad Y_i|\theta_i\sim \Poi(\theta_i),i\in [n].
	$$
	In view of the standard inequality that supremum over $\theta_1,\dots\theta_n$ is greater than the expected value with respect to $\theta_i\sim \hat G,i\in[n]$, we get
	$$
	\sup_{\theta_1,\dots,\theta_n} \sum_{i=1}^n
	\qth{\pth{\theta_i-\hat \theta_{\hat G}(Y_i)}^2}
	= \EE_{\theta_i\sim \hat G,i\in[n]}\sum_{i=1}^n\qth{\pth{\theta_i-\hat \theta_{\hat G}(Y_i)}^2},
	\quad Y_i|\theta_i\sim \Poi(\theta_i),i\in [n].
	$$
	As $\sth{\hat\theta_{\hat G}(Y_1),\dots,\hat\theta_{\hat G}(Y_n)}$ is the Bayes estimator with respect to the prior $\theta_i\inddistr\hat G,i\in n$, given the data generating model $Y_i|\theta\sim\Poi(\theta_i),i\in [n]$, we use \cite[Section 5.1, Theorem 1.4]{lehmann2006theory} to conclude that $\sth{\hat\theta_{\hat G}(Y_1),\dots,\hat\theta_{\hat G}(Y_n)}$ is a minimax estimator.
	\end{proof}}
	
	\section{Proof of \prettyref{thm:existence-uniqueness}}
	
	\label{app:EB_existence-uniqueness}
	
	We first prove the result for the constrained solution $\argmin_{Q\in\calP([0,h])} \dist(p\|f_{Q})$. 
	As mentioned towards the end of the proof, this also implies the desired result for the unconstrained solution.
	Suppose that $p$ is supported on $\{y_1,\dots,y_m\} \subset \integers_+$. 
	Define
	\begin{align}\label{eq:EB_mS}
		S\triangleq\sth{(f_{Q}(y_1),\dots,f_Q(y_m)): Q\in \calP([0,h])},
	\end{align} 
	where $f_Q(y) = \Expect_{\theta \sim Q}[f_\theta(y)]$ is the probability mass function of the Poisson mixture \prettyref{eq:EB_mixture}, and 
	$f_\theta(y)=e^{-\theta} \theta^y/y!$.
	We claim that $S$ is convex and compact.\footnote{In this case, $S$ is in fact the closed convex hull of the set $\{(f_{\theta}(y_1),\dots,f_\theta(y_m)): \theta\in [0,h]\}$.}
	The convexity follows from definition. For compactness, note that $S$ is bounded since $\sup_{\theta\geq 0} f_\theta(y) = e^{-y} y^y/y!$, so it suffices to check $S$ is closed. Let 
	$(f_1',\ldots,f_m')\in\reals_+^m$ be the limiting point of $(f_{Q_k}(y_1),\dots,f_{Q_k}(y_m))$ for some sequence 
	$\{Q_k\}$ in $\calP([0,h])$. 
	By Prokhorov's theorem, there is a subsequence $\{Q_{k_\ell}\}$ that converges weakly to some $Q'\in \calP([0,h])$.
	Since $\theta \mapsto f_\theta(y)$ is continuous and bounded, we have $f_j'=f_{Q'}(y_j)$ for all $j$. In other words, $S$ is closed.

	Next, define $v:S \to \reals$ by $v(f_1,\dots,f_m)=\sum_{i=1}^m
	{\ell(p(y_i),f_i)}$.
	By \prettyref{pt:EB_t-ell}, the value of the min-distance optimization can be written as
	\begin{equation}
		\min_{Q\in\calP([0,h])} \dist(p\|f_{Q}) = t(p) + \min_{(f_1,\ldots,f_m) \in S} v(f_1,\ldots,f_m).
		\label{eq:mindist-equiv}
	\end{equation}
	Furthermore, by assumption $\ell(0,b)\equiv 0$ and $b\mapsto \ell(a,b)$ is strictly convex for $a>0$. Thus $v$ is strictly convex.
	Therefore, there exists a unique point $(f_1^*,\dots,f_m^*) \in S$ that achieves the minimum on the right side of \prettyref{eq:mindist-equiv}. 
	Thus, 
	the left side has a minimizer $\hatGn\in\calP([0,h])$ that satisfies 	
	\begin{align}\label{eq:EB_m4}
		f_{\hatGn}(y_j)=f_j^*,j=1,\dots,m
	\end{align}
	It remains to show that the above representation is unique at the special point $(f_1^*,\dots,f_m^*)$; this argument relies on the specific form of the Poisson density.
	Let ${\hatGn}$ be one such minimizer.
	By the first-order optimality condition (see \prettyref{eq:EB_first-order-optimality} in \prettyref{sec:first-order-optimality}),
	\begin{align}\label{eq:opt1}
		D_{\hatGn}(\theta)&=\sum_{i=1}^m a_i
		(f_\theta(y_i)-f_i^*)\geq 0, \quad \forall 0\leq\theta\leq h; 
		\nonumber\\
		D_{\hatGn}(\theta)&=0, \quad \text{for $\hatGn$-almost every $\theta$},
	\end{align}
	where $a_i \triangleq \frac d{df}\ell(p(y_i),f)|_{f=f_i^*} < 0$, since $\ell$ is strictly decreasing in the second coordinate and $f_i^*>0$.
	Define 
	\begin{align*}
		b_i={a_i \over \sum_{i=1}^m a_i f_i^*} > 0.
	\end{align*}
	
	As $\ell$ is strictly decreasing in second coordinate, $\frac d{df}\ell(p(y_i),f)< 0$ for all $f\in\reals_+,i=1,\dots,m$. Using this, we rearrange  \eqref{eq:opt1} to get
	\begin{align}\label{eq:EB_support-optimality}
		\sum_{i=1}^m {b_i\over y_i!}\theta^{y_i}&\leq  e^{\theta},\forall \theta\in [0,h] 
		, \nonumber\\
		\sum_{i=1}^m {b_i\over y_i!}\theta^{y_i}&= e^{\theta}
		\text{ for each $\theta$ in the support of } \hatGn.
	\end{align}
	Then the following lemma shows that the support of $\hatGn$ has at most $m$ points.
	\begin{lemma}
		\label{lmm:EB_roots}
		Suppose that $\sum_{i=1}^m \beta_i\theta^{y_i}\leq e^{\theta}$ for all $\theta\in[0,h]$ where $\beta_i \in\reals$ and $h>0$. Then 
		the number of solutions to $\sum_{i=1}^m \beta_i\theta^{y_i}= e^{\theta}$ in $\theta \in [0,h]$ is at most $m$.
	\end{lemma}
	\begin{proof}
		The proof is a modification of \cite[Lemma 3.1(2)]{S76}, which deals with the specific case $h=\infty$. 
		Recall the following version of Descartes' rule of signs \cite[Part V, Problem 38 and 40]{polya1998problems}: 		
		Consider an entire function (i.e., a power series whose radius of convergence is infinity) $\phi(x)=a_0+a_1x+a_2x^2+\dots$ with real coefficients. Let $r$ be the number of strictly positive zeros of $\phi$ counted with their multiplicities and let $s$ be the number of sign changes\footnote{The number of sign changes is the number of pairs $0\leq i<j$ such that $a_ia_j<0$ and either $j=i+1$ or $a_k=0$ for all $i<k<j$.} in the sequence $a_0,a_1,\dots$. Then $r\leq s$. 
		We apply this fact to the function
		$$\phi(\theta)
		=\sum_{i=1}^m \beta_i \theta^{y_i}- e^{\theta}
		=\sum_{j=0}^\infty a_j\theta^j,$$
		where
		\[
		a_j=
		\begin{cases}
			\beta_{i}-\frac 1{y_i!}		& j =y_i, i=1,\dots, m \\
			-\frac 1{j!}<0	 & \text{else}\\
		\end{cases}
		\]
		
		\paragraph{Case 1:} Suppose that 0 is a root of $\phi(\cdot)$. Then $a_0=0$. As there are at most $m-1$ positive coefficients in $a_0,a_1,\dots$, there can be at most $2(m-1)$ sign changes, which implies at most $2(m-1)$ positive roots of $s$ counting multiplicities. Note that, as $\phi(\theta)\indc{\theta\in (0,h)}\leq 0$ 
		and $s$ is an entire function, each root of $s$ inside $(0,h)$ has multiplicity at least 2. Suppose that $m_h$ is the multiplicity of $h$ as a root of $\phi(\cdot)$, which we define to be 0 when $h$ is not a root. This means that the total number of distinct roots in $(0,h)$ is at most the largest integer before $(2(m-1)-m_h)/2$. If $h$ is not a root, then the number of distinct roots in $(0,h)$ is at most $m-1$. If $h$ is a root, then its multiplicity is at least 1, and hence, the number of distinct roots in $(0,h)$ is at most $m-2$. Hence, there are at most $m$ many distinct roots in $[0,h]$.  
		
		\paragraph{Case 2:} Suppose that 0 is not a root of $\phi(\cdot)$. As there are at most $m$ positive coefficients in $a_0,a_1,\dots$, there can be at most $2m$ sign changes, which implies at most $2m$ positive roots counting multiplicities. By a similar argument as in the previous case, the total number of distinct roots in $(0,h)$ is at most the largest integer before $(2m-m_h)/2$. If $h$ is not a root, then the number of distinct roots in $(0,h)$ is at most $m$. If $h$ is a root, then the number of distinct roots in $(0,h)$ is at most $m-1$. Hence, in total, there are at most $m$ distinct roots in $[0,h]$.  
		\end{proof}
		Suppose that there are $r(\leq m)$ different $\theta_i$'s (denote them by $\theta_1\dots,\theta_r$) for which \eqref{eq:EB_support-optimality} holds. This implies given any optimizer ${\hatGn}$ its atoms form a subset of $\sth{\theta_1\dots,\theta_r}$. Let $w_j$ be the weight ${\hatGn}$ puts on $\theta_j$. Then in view of \eqref{eq:EB_m4} we get that 
		$$
		\sum_{j=1}^r{w_je^{-\theta_j}}\theta_j^{y_i}=f_i^*y_i!,\quad i=1,\dots,r.
		$$
		The matrix $\{\theta_j^{y_i}:j=1,\dots,r,i=1,\dots,m\}$ has full column rank, and hence the vector $(w_1,\dots,w_r)$ can be solve uniquely. This also implies the uniqueness of the optimizer ${\hatGn}$. This finishes the proof for the constrained solution.
		
		Next we argue for the unconstrained minimizer $\argmin_{Q}\dist(p\|f_Q)$. In view of \prettyref{lmm:EB_support-g-hat} below, we get that the unconstrained minimum-distance estimator is supported on $[0,h]$ with $h=\max_{i=1,\dots,m}y_i$. Then, from the above proof for $\argmin_{Q\in\calP([0,h])}\dist(p\|f_Q)$, the existence and uniqueness of the unconstrained estimator follow.


	
	\begin{lemma}\label{lmm:EB_support-g-hat}
		Let $\dist$ satisfy \prettyref{pt:EB_t-ell} and let $p$ be a probability distribution on $\integers_+$ with support $\{y_1,\dots,y_m\}$. Then the minimizer $\argmin_Q \dist(p\|f_Q)$ is supported on the interval $[y_{\min},y_{\max}]$, where 
		$y_{\min}=\min_{i=1,\dots,m}y_i,y_{\max}=\max_{i=1,\dots,m}y_i$. 
	\end{lemma}
	
	\begin{proof}
		
		Let $Q$ be a distribution with $Q([0,y_{\min}))+Q((y_{\max},\infty))>0$. Define another distribution $\tilde Q$ by
		$$\tilde Q(\cdot)=
		Q([0,y_{\min})) \delta_{y_{\min}}(\cdot)
		+Q(\cdot\cap [y_{\min},y_{\max}])
		+Q((y_{\max},\infty)) \delta_{y_{\max}}(\cdot).
		$$
		In other words, $\tilde Q$ moves the masses of $Q$ on the intervals $[0,y_{\min})$ (resp.~$(y_{\max},\infty)$) to the point $y_{\min}$ (resp.~$y_{\max}$). As $f_\theta(y)$ is strictly increasing in $\theta\in[0,y)$ and strictly decreasing in $\theta\in (y,\infty)$ we get for each $i=1,\dots,n$
		\begin{align*}
			f_Q(y_i)
			&=\int f_{\theta}(y_i)dQ(\theta)
			\nonumber\\
			&=\int_{0\leq \theta< y_{\min}} f_{\theta}(y_i)dQ(\theta)
			+\int_{y_{\min}\leq \theta\leq y_{\max}} f_{\theta}(y_i)dQ(\theta)
			+\int_{y_{\max}>\theta} f_{\theta}(y_i)dQ(\theta)
			\nonumber\\
			&<
			Q([0,y_{\min}))
			f_{y_{\min}}(y_i)
			+\int_{y_{\min}\leq\theta\leq 	y_{\max}} f_{\theta}(y_i)dQ(\theta)
			+Q((y_{\max},\infty))
			f_{y_{\max}}(y_i)
			\nonumber\\
			&
			=\int f_{\theta}(y_i)d\tilde Q(\theta)
			=f_{\tilde Q}(y_i).
		\end{align*}
		Hence, by \prettyref{pt:EB_t-ell}, we get
		\begin{align}
			&\dist(p\|f_Q)
			=t(p)+ \sum_{y\geq 0}
			\ell(p(y),f_Q(y))
			\nonumber\\
			&~\stepa{=}t(p)+ \sum_{y:p(y)>0}
			{\ell(p(y),f_Q(y))}
			\stepb{>}
			t(p)+
			\sum_{y:p(y)>0}{\ell\pth{p(y),f_{\tilde Q}(y)}}
			=\dist(p\|f_{\tilde Q}),
			\label{eq:EB_ell-representation}
		\end{align}
		where (a) follows from $\ell(0,\cdot)=0$;
		(b) follows as the function $b\mapsto \ell(a,b)$ is strictly decreasing.  In other words, given any $Q$ with $Q([0,y_{\min}))+Q((y_{\max},\infty))>0$ we can produce $\tilde Q$ supported on $[y_{\min},y_{\max}]$ such that $\dist(p\|f_{\tilde Q})< \dist(p\|f_Q)$. Hence, the claim follows.
	\end{proof}

	\section{Proof of \prettyref{lmm:EB_general-regret-bound}}
	\label{app:EB_general-regret-bound}
	
	Let $\theta\sim G, Y|\theta\sim f_{\theta}$. Then for any $\hatGn$ independent of $Y$, we can write
	$\Regret(\hat G;G)=\sum_{y=0}^\infty (\hat\theta_{\hat G}(y)-\hat \theta_{G}(y))^2f_G(y)=\EE_G\qth{\pth{\hat\theta_{{\hatGn}}(Y)-\hat\theta_{ G}(Y)}^2}$; cf.~\prettyref{eq:regret}.
	Fix $h>0$ and note the following
	\begin{itemize}
		\item 
		~$
			\text{mmse}(G)
			=\EE[(\hat\theta_G - \theta)^2]
			\geq \PP\qth{\theta\in [0,h]} \ \EE_{\theta\sim G}[(\hat\theta_G - \theta)^2 | \theta\in [0,h]]
			\ge \PP\qth{\theta\in [0,h]} \mmse(G_h)
		$
		\item ~$\text{mmse}(G)\leq \sqrt{\EE_G[\theta^4]}\leq \sqrt{M}$, and
		\item ~For any fixed distribution $\hat G$
		\begin{align}
			\EE_{G}\qth{(\hat\theta_{{{\hat G}}}(Y)-\theta)^2}
			&\leq 
			\EE_{G}\qth{(\hat\theta_{{{\hat G}}}(Y)-\theta)^2\indc{\theta\leq h}}
			+\EE_{G}\qth{(\hat\theta_{{{\hat G}}}(Y)-\theta)^2\indc{\theta>h}}
			\nonumber\\
			&\stepa{\leq}
			\EE_{G}\qth{\left.(\hat\theta_{{{\hat G}}}(Y)-\theta)^2\right|{\theta\leq h}}
			+\sqrt{\EE_{G}\qth{(\hat\theta_{{{\hat G}}}(Y)-\theta)^4}\EE_G\qth{\indc{\theta>h}}}
			\nonumber\\
			&\stepb{\leq} \EE_{G_{h}}\qth{(\hat\theta_{{{\hat G}}}(Y)-\theta)^2}+\sqrt{8(\hat h^4+\EE_{G}[\theta^4])G((h,\infty))}
			\nonumber\\
			&=\EE_{G_{h}}\qth{(\hat\theta_{{{\hat G}}}(Y)-\theta)^2}
			+\sqrt{8(\hat h^4+M)G((h,\infty))}.
			\label{eq:EB_m7}
		\end{align}
		where step (a) followed by Cauchy-Schwarz inequality and step (b) followed as $(x+y)^4\leq 8(x^4+y^4)$ for any $x,y\in \reals$.
	\end{itemize}
	Using these we get
	\begin{align}\label{eq:EB_m1}
		\Regret(\hatGn;G)
		&=
		\EE_{G}\qth{(\hat\theta_{{{\hat G}}}(Y)-\theta)^2}-\text{mmse}(G)
		\nonumber\\
		&\leq 
		\EE_{G_{h}}\qth{(\hat\theta_{{{\hat G}}}(Y)-\theta)^2}-\text{mmse}(G_{h})
		+\text{mmse}({G_{h}})-\text{mmse}(G)
		+\sqrt{8(\hat h^4+M)G((h,\infty))}
		\nonumber\\
		&\leq \EE_{G_{h}}\qth{(\hat\theta_{{{\hat G}}}(Y)-\hat\theta_{G_h}(Y))^2}
		+\pth{\frac 1{G([0,h])}-1}{\mmse(G)}
		+\sqrt{8(\hat h^4+M)G((h,\infty))}
		\nonumber\\
		&\leq \EE_{G_{h}}\qth{(\hat\theta_{{{\hat G}}}(Y)-\hat\theta_{G_h}(Y))^2}
		+{\frac {G((h,\infty))}{G([0,h])}}\sqrt M
		+\sqrt{8(\hat h^4+M)G((h,\infty))}
		\nonumber\\
		&\leq \EE_{G_{h}}\qth{(\hat\theta_{{{\hat G}}}(Y)-\hat\theta_{G_h}(Y))^2}
		+{(1+2\sqrt 2)\sqrt{(\hat h^4+M)G((h,\infty))}\over G([0,h])}.
	\end{align}
	Next we bound the first term. Fix $K\geq 1$. Using $\hat\theta_{{G_h}}(y)\leq h,\hat\theta_{\hat {G}}(y)\leq \hat h$ we have
	\begin{align*}
		&\EE_{G_h}\qth{(\hat\theta_{\hat{G}}(Y)-\hat\theta_{{G}_h}(Y))^2\indc{Y\leq K-1}}
		\nonumber\\
		&=\sum_{y=0}^{K-1}(y+1)^2f_{G_h}(y)\pth{{f_{\hat{G}}(y+1)\over f_{\hat{G}}(y) }-{f_{G_h}(y+1)\over f_{G_h}(y)}}^2
		\nonumber\\
		&\stepa{\leq} \sum_{y=0}^{K-1}(y+1)^2f_{G_h}(y)
		\left\{3\pth{{f_{\hat{G}}(y+1)\over f_{\hat{G}}(y) }-{2f_{\hat{G}}(y+1)\over f_{G_h}(y)+f_{\hat{G}}(y)}}^2+3\pth{{f_{G_h}(y+1)\over f_{G_h}(y) }-{2f_{G_h}(y+1)\over f_{G_h}(y)+f_{\hat{G}}(y)}}^2
		\right.
		\nonumber\\
		&\quad \left.+3\pth{2f_{G_h}(y+1)-2f_{\hat{G}}(y+1)\over f_{G_h}(y)+f_{\hat{G}}(y)}^2\right\}
		\nonumber\\
		&\leq 3\sum_{y=0}^{K-1}\left\{
		\pth{(y+1)f_{\hat{G}}(y+1)\over f_{\hat{G}}(y) }^2{(f_{G_h}(y) -f_{\hat{G}}(y) )^2\over f_{G_h}(y)+f_{\hat{G}}(y)}
		+\pth{(y+1)f_{G_h}(y+1)\over f_{G_h}(y)}^2{(f_{G_h}(y) -f_{\hat{G}}(y) )^2\over f_{G_h}(y)+f_{\hat{G}}(y)}\right.
		\nonumber\\
		&\quad \left.+4(y+1)^2{(f_{G_h}(y+1)-f_{\hat{G}}(y+1))^2
			\over f_{G_h}(y)+f_{\hat{G}}(y)}\right\}
		\nonumber\\
		&= 3(\{\hat\theta_{G_h}(y)\}^2+\{\hat\theta_{\hat G}(y)\}^2)\sum_{y=0}^{K-1}{(f_{G_h}(y)-f_{\hat{G}}(y))^2\over f_{G_h}(y)+f_{\hat{G}}(y)}
		+12\sum_{y=0}^{K-1}(y+1)^2{(f_{G_h}(y+1)-f_{\hat{G}}(y+1))^2
			\over f_{G_h}(y)+f_{\hat{G}}(y)}
		\nonumber\\
		&\leq 3(h^2+\hat h^2)\sum_{y=0}^{K-1}{(f_{G_h}(y)-f_{\hat{G}}(y))^2\over f_{G_h}(y)+f_{\hat{G}}(y)}
		+12\sum_{y=0}^{K-1}(y+1)^2{(f_{G_h}(y+1)-f_{\hat{G}}(y+1))^2
			\over f_{G_h}(y)+f_{\hat{G}}(y)}
	\end{align*}
	where (a) followed from $(x+y+z)^2\leq 3(x^2+y^2+z^2)$ for any $x,y,z\in\reals$. Using
	$(\sqrt {f_{G_h}(x)}+\sqrt {f_{\hat{G}}(x)})^2\leq 2(f_{G_h}(x)+f_{\hat{G}}(x))$ for $x=y,y+1$ we continue the last display to get
	\begin{align*}
		&\EE_{G_h}\qth{(\hat\theta_{{\hat G}}(Y)-\hat\theta_{G_h}(Y))^2\indc{Y\leq K-1}}
		\nonumber\\
		&\leq 6(h^2+\hat h^2)\sum_{y=0}^{K-1}(\sqrt{f_{G_h}(y)}-\sqrt{f_{\hat{G}}(y)})^2
		\nonumber\\
		&\quad +24K\max_{y=0}^{K-1}{(y+1)f_{G_h}(y+1)+(y+1)f_{\hat{G}}(y+1)\over f_{G_h}(y)+f_{\hat{G}}(y)}
		\sum_{y=0}^{K-1}(\sqrt{f_{G_h}(y+1)}-\sqrt{f_{\hat{G}}(y+1)})^2
		\nonumber\\
		&\leq \pth{6(h^2+\hat h^2)+24(h+\hat h)K} {H^2(f_{\hat{G}},f_{G_h})}.
	\end{align*}
	Again using $\hat\theta_{{G_h}}(y)\leq h,\hat\theta_{\hat {G}}(y)\leq \hat h$ we bound $\EE_{G_h}\qth{(\hat\theta_{{\hat G}}(Y)-\hat\theta_{G_h}(Y))^2\indc{Y\geq K}}$ by $(h+\hat h)^2\epsilon_K(G_h)$. Combining this with the last display we get
	\begin{align*}
		\EE_{G_h}\qth{(\hat\theta_{{\hat G}}(Y)-\hat\theta_{G_h}(Y))^2}\leq \sth{6(h^2+\hat h^2)
			+24(h+\hat h)K} {H^2(f_{\hat{G}},f_{G_h})}+(h+\hat h)^2\epsilon_K(G_h).
	\end{align*}
	In view of above continuing \eqref{eq:EB_m1}
	we have
	\begin{align}
		\Regret(\hat G;G)
		\leq  
		&\sth{6(h^2+\hat h^2)
			+24(h+\hat h)K}{H^2(f_{\hat{G}},f_{G_h})}
		\nonumber\\
		&\quad +(h+\hat h)^2\epsilon_K(G_h)
		+{(1+2\sqrt 2)\sqrt{(M+\hat h^4)G((h,\infty))}\over G([0,h])}.
		\label{eq:EB_m9}
	\end{align}
	Using triangle inequality and $(x+y)^2\leq 2(x^2+y^2)$ we get
	\begin{align}
		H^2(f_{{\hatGn}},f_{G_h})
		\leq 2\sth{H^2(f_{G},f_{{\hatGn}})+H^2(f_{G_h},f_{G})}.
		\label{eq:EB_m10}
	\end{align}
	Note that
	\[
	H^2(f_{G_h},f_{G}) \leq 2 \TV(f_{G_h},f_{G}) \leq 2 \TV(G_h,G) = 4 G((h,\infty)).
	\]
	where $\TV$ denotes the total variation and the middle inequality applies the data-processing inequality \cite{Csiszar67} and the last equality followed as
	\begin{align*}
		\TV(G_h,G)
		&=\int_0^h\abs{dG_h(\theta)-dG(\theta)}+\int_h^\infty
		dG(\theta)
		\nonumber\\
		&=\pth{\frac 1{G([0,h])}-1}\int_0^h dG(\theta)+G((h,\infty))
		=2G((h,\infty)).
	\end{align*}
	
	Then, combining \eqref{eq:EB_m9}, \eqref{eq:EB_m10} and using $\epsilon_K(G_h)\leq {\epsilon_K(G)\over G([0,h])}$ we get the desired bound
	\begin{align*}
		\Regret(\hat G;G)
		&\leq \sth{12(h^2+\hat h^2)
			+48(h+\hat h)K}\pth{H^2(f_{\hat{G}},f_{G})+4G((h,\infty))}
		\nonumber\\
		&\quad +(h+\hat h)^2{\epsilon_K(G)\over G([0,h])}
		+{(1+2\sqrt 2)\sqrt{(M+\hat h^4)G((h,\infty))}\over G([0,h])}.
	\end{align*}

	\section{Auxiliary results}
	\label{app:EB_subexpo-properties}
	
	\begin{lemma}\label{lmm:anru-zhang-results-poisson}
		Let $Y\sim \Poi(h)$. Then we have
		\begin{align*}
			\PP[Y>h+y]\leq e^{-{y^2\over 4h}}
			~\forall~ 0<y<\frac h2,\quad 
			\PP[Y<h-y]\leq e^{-{5y^2\over 9h}} ~\forall~ 0<y<\frac h3
		\end{align*}
	\end{lemma}	
	\begin{proof}
		~The proof of the above result follows from  \cite[Section]{zhang2020non} with the inequality $\log(1+t)\geq (t-2t^2/3)\indc{t\in (-\frac 13,0)}+(t-t^2/2)\indc{t\in(0,\frac 12)}$.
	\end{proof}
	
	\begin{lemma}\label{lmm:bounded_prior_properties}
		Let $h>0$ and $G\in\calP([0,h])$. If $\sth{Y_i}_{i=1}^n\iiddistr f_G$, then the following are satisfied given  any $t\geq 0,a\geq 1, n\geq 3$ and $K=\min\sth{{a(he^2+2)\log n\over \log\log n},he^2+a\log n}$.
		$$ \PP\qth{Y_1>K+t}\leq {2e^{-t}\over n^a},\quad \EE[Y_{\max}^4]\leq 5K^4+ 8.
		$$
	\end{lemma}

	\begin{proof}
		~Let $G\in \calP([0,h])$. 
		As $p(\theta)=e^{-\theta}\theta^y$ is increasing in $\theta\in [0,y],y>0$, for $\ell>2h$,
		\begin{align}\label{eq:r3}
			\PP[Y_1\geq \ell]
			=\sum_{y=\ell}^\infty
			\int_0^h{e^{-\theta}\theta^y\over y!}
			G(d\theta)
			\leq \sum_{y=\ell}^\infty
			{e^{-h}h^y\over y!}
			\leq {h^\ell\over \ell!}
			{\sum_{y-\ell=0}^\infty\pth{\frac h\ell}^{y-\ell}}
			\leq 2\pth{he\over \ell}^{\ell}.
		\end{align}
		where the last inequality followed using $\ell!\geq \pth{\frac {\ell}e}^{\ell}$ from the Stirling's formula. To get to our results, first, let $K_1={a(he^2+2)\log n\over \log\log n}, a\geq 1$. Using the fact $\log\log\log n<\frac {\log\log n}2$ for all $n\geq 3$ we continue the last display with $\ell=K_1+t$ to get
		\begin{align}
			\PP[Y_1\geq K_1+t]
			&\leq 2\pth{he\over K_1+t}^{K_1}\pth{he\over K_1+t}^{t}
			\leq 2\pth{he\over K_1}^{K_1} e^{-t} 
			\nonumber\\
			&\leq 2e^{-t}\pth{\log\log n\over \log n}^{a(he^2+2)\log n\over \log\log n}
			\nonumber\\
			&\leq 2e^{-t-\pth{\log\log n-\log\log\log n}{a(he^2+2)\log n\over \log\log n}}
			\leq 2e^{-t-a\log n}
			\leq \frac {2e^{-t}}{n^a}
			\label{eq:EB_epK_1},
		\end{align}
		as required. Next, considering $K_2=he^2+a\log n$ we continue \eqref{eq:r3} with $\ell=K_2+t$ to get
		\begin{align}
			\PP[Y_1\geq K_2+t]
			&\leq 2\pth{he\over K_2}^{K_2}\pth{he\over K_2+t}^{t}
			\leq 2\pth{1-{a\log n\over K_2}}^{K_2}e^{-t}
			\leq 2e^{-a\log n} e^{-t}
			\leq \frac {2e^{-t}}{n^a}
			\label{eq:EB_epK_2}.
		\end{align}
		Choosing $K=\min\{K_1,K_2\}$ we get the desired result.
		
		Next we bound $\EE[Y_{\max}^4]$. For any nonnegative integer valued random variable  $Z$, using 
		\begin{align}\label{eq:E-z4}
			\EE[Z^4]=\sum_{z\geq 1} z^4\PP[Z=z]
		\leq 4\sum_{z\geq 1} \sum_{k=1}^z k^3\PP[Z=z]
		=4\sum_{k\geq 1}k^3\sum_{z\geq k}\PP[Z=z]
		=4\sum_{k\geq 1}k^3\PP[Z\geq k] 
		\end{align}
		we have
		\begin{align*}
			\EE\qth{\pth{Y_{\max}}^4}
			&\leq 4\sum_{y} y^3\PP
			\qth{Y_{\max}>y}
			\nonumber\\
			&\leq 4K^4
			+n\sum_{y\geq K+1} y^3 \PP\qth{Y_1>y}
			\nonumber\\
			&= 4K^4
			+n\sum_{t\geq 1} (K+t)^3 \PP[Y_1>K+t]
			\nonumber\\
			&\stepa{\leq} 4K^4
			+{4\over n^4}\sum_{t\geq 1}(K^3+t^3)e^{-t}
			\leq 4K^4
			+{\frac {K^3}4}\sum_{t\geq 1}e^{-t}
			+\sum_{t\geq 1}t^3e^{-t}
			\leq 5K^4+ 8,
		\end{align*}
		 where (a) followed from $(x+y)^3\leq 4(x^3+y^3), x,y\geq 0$.
	\end{proof}
	
	\begin{lemma}\label{lmm:subexpo_properties}
		Given any $s>0$ and $G\in\subexpo(s)$, the following are satisfied.
	\begin{enumerate}[label=(\roman*)~]
		\item \label{pt:i} If $\theta\sim {G}$, then $\EE[\theta^4]\leq 12s^4$.
		\item \label{pt:ii} If $\sth{Y_i}_{i=1}^n\iiddistr f_G$, then for any $K\geq 1$
		\begin{align}
			\label{eq:EB_prop_mixture}
			\PP\qth{Y_1\geq K}
			\leq 3 e^{-K\log\pth{1+\frac 1{2s}}},
			\quad
			\EE\qth{Y_{\max}^4}\leq {64(\log n)^4+90\over \pth{\log\pth{1+\frac 1{2s}}}^4}.
		\end{align}
	\end{enumerate}
	\end{lemma}
	\begin{proof}
		
	To prove (i) we note that for any $s,M>0$ using integral by parts we have $\int_{0}^Mx^3e^{-\frac xs} dx= \qth{-se^{-\frac xs}(6s^3+6s^2x+3sx^2+x^3)}_{0}^M$. Then we get using the definition of $\subexpo(s)$ tail probabilities and with limit as $M\to \infty$
	$$
	\EE[\theta^4]
	=4\int_0^\infty y^3\PP
	\qth{\theta>y}dy
	\leq 2\int_0^\infty y^3 e^{-\frac ys}dy
	= 2 \lim_{M\to \infty}\qth{-se^{-\frac ys}(6s^3+6s^2y+3sy^2+y^3)}_{0}^M
	\leq 12s^4.
	$$
	The proof of the property (ii) is as follows. Using $\EE_{Z\sim\Poi(\theta)}\qth{e^{Zt}}=e^{\theta(e^t-1)},t>0$ and denoting $c(s)=\log{1+2s\over 2s}$ we have 
	\begin{align*}
		&\EE\qth{e^{Y_1c(s)}}
		=\EE_{\theta\sim {G}}\qth{\EE_{Y_1\sim \Poi(\theta)}\qth{\left.e^{Y_1c(s)}\right|\theta}}
		=\EE_{G}\qth{e^{\theta\over 2s}}
		=\int_0^\infty e^{\theta/2s}\ G(d\theta)
		\nonumber\\
		&=1+\int_{\theta=0}^\infty \int_{x=0}^\theta {e^{x/2s}\over 2s}dx\ G(d\theta)
		=1+\int_{x>0} {e^{x/2s}\over 2s}\ G([x,\infty))dx
		\stepa{\leq} 
		1+\int_{x>0}{e^{-x/2s}\over s}dx
		\leq 3
	\end{align*}
	where (a) followed by using tail bound for $\subexpo(s)$ distribution $G$. In view of Markov inequality 
	\begin{align*}
		\PP\qth{Y_1\geq K}
		\leq \EE\qth{e^{Y_1c(s)}}e^{-c(s)K}
		\leq
		3e^{-K\log\pth{1+\frac 1{2s}}}.
	\end{align*}
	The expectation term is bounded as below. Pick $L$ large enough such that $v(y) = y^3e^{-y\log\pth{1+\frac 1{2s}}}$ is decreasing for all $y\geq L$. Then we can bound $\sum_{y=L+1}^{\infty} v(y)\leq \int_{y>L} v(y) dy$. Then, using the last inequality for such $L>0$ and \eqref{eq:E-z4}
	\begin{align*}
		\EE\qth{\pth{Y_{\max}}^4}
		&\leq 4\sum_{y} y^3\PP
		\qth{Y_{\max}>y}
		\nonumber\\
		&\leq 4L^4
		+n\sum_{y\geq L+1} y^3 \PP\qth{Y_1>y}
		\nonumber\\
		&{\leq} 4L^4
		+{3n}\int_{y\geq L} y^3e^{-y\log\pth{1+\frac 1{2s}}} dy
		\nonumber\\
		&\stepa{\leq} 4L^4+\frac {3n}{\sth{\log\pth{1+\frac 1{2s}}}^4}\int_{z>L\log\pth{1+\frac 1{2s}}} z^3e^{-z}dz
		\nonumber\\
		&\stepb{\leq} 4L^4+\frac {45n}{\sth{\log\pth{1+\frac 1{2s}}}^4}\int_{z>L\log\pth{1+\frac 1{2s}}} e^{-z/2}dz
		\leq 4L^4+{90ne^{-\frac L2\log\pth{1+\frac 1{2s}}}\over \sth{\log\pth{1+\frac 1{2s}}}^4},
	\end{align*} 
	where (a) followed from a change of variable, (b) followed using $x^3\leq 15e^{-\frac x2}$ for any $x>0$. Choosing $L={2\log n\over \log\pth{1+\frac 1{2s}}}$ we get the desired result.
\end{proof}
	
	\section{Proofs of the multidimensional results}
	
	\label{app:multidim-proofs}
	
	\subsection{Density estimation in multiple dimensions}
	
	We will assume $d\leq n$, as the results are vacuous otherwise. The proof of \prettyref{thm:density_multidim} is based on a similar truncation idea as in the proof of \prettyref{thm:EB_density-H2}. For this section, we will use the notation
	\begin{align}
		\label{eq:K-defn}
		K=\begin{cases}
			\min\sth{\max\{1,h\}{\log n\over \log\log n},h+\log n},& {G\in \calP([0,h]^d},\\
			\max\{1,s\}{\log n}, & {\text{marginals of $G$ are in $\subexpo(s)$}},
		\end{cases}
	\end{align}
	unless we specify it differently. We first note the following result.
	
	\begin{lemma}
		\label{lmm:support-multidim}
		There exists absolute constant $\tilde c_1$ such that the following holds. If $Y\sim f_G$ such that $G\in \calP([0,h]^d)$ or all the marginals of $G$ belong to $\subexpo(s)$, then $\PP\qth{\vY\notin [0,\tilde c_1 K]^d}\leq \frac d{n^{10}}$.
	\end{lemma}
	\begin{proof}
		
		 ~From the proof of \prettyref{lmm:bounded_prior_properties} and \prettyref{lmm:subexpo_properties}, we get that there exists a constant $\tilde c_1$ such that with probability at least $1-\frac 1{n^{10}}$ all the coordinates of the random variable $\vY$ lie within $[0,\tilde c_1 K]$. Then using a union bound over all the coordinates we achieve the desired result. 
	\end{proof}

	\begin{proof}[Proof of \prettyref{thm:density_multidim}~]
		Suppose that the $\dist$ function, for which we compute the minimum distance estimator, satisfies \prettyref{pt:EB_multi_sandwich}. Then, using a proof strategy identical to proving the result \eqref{eq:EB_m6} in the one-dimensional case, we get
		\begin{align*}
			H^2(f_{G},f_{{\hatGn}})
			\leq \frac 2{c_1}{(\dist(p^\sfem_n\|f_{{\hatGn}})+\dist(p^\sfem_n\|f_{G}))}
			\leq \frac 4{c_1}{\dist(p^\sfem_n\|f_{G})}.
		\end{align*}
		In view of \prettyref{pt:EB_multi_sandwich} we bound ${1\over c_2}\dist$ by $\chi^2$ and use $R=\tilde c_1 K$ to get the following
		\begin{align}
			&~\frac{1}{c_2} \EE\qth{\dist(p^\sfem_n\|f_{G})\indc{\vY_i \in [0,R]^d\ \forall i=1,\dots,n}}
			\nonumber\\
			&\leq \EE\qth{\chi^2(p^\sfem_n\|f_{G})\indc{\vY_i \in [0,R]^d\ \forall i=1,\dots,n}} 
			=
			\sum_{\vy}{\EE\qth{(p^\sfem_n(\vy)-f_{G}(\vy))^2\indc{\vY_i \in [0,R]^d\ \forall i=1,\dots,n}}\over f_{G}(\vy)}
			\nonumber\\
			&\stepa{=} 
			\sum_{\vy \in [0,R]^d }{\EE\qth{(p^\sfem_n(\vy)-f_{G}(\vy))^2\indc{\vY_i \in [0,R]^d\ \forall i=1,\dots,n}}\over
				f_{G}(\vy)}  +\sum_{\vy\notin [0,R]^d}f_{G}(\vy) \PP[\vY_i \in [0,R]^d\ \forall i=1,\dots,n]
			\nonumber\\
			&\leq
			\sum_{\vy \in [0,R]^d }{\EE\qth{(p^\sfem_n(\vy)-f_{G}(\vy))^2}\over
				f_{G}(\vy)}  + \PP_{\vY\sim f_G}[\vY \notin [0,R]^d]
			\nonumber\\
			&\stepb{\leq}
			\frac{1}{n}\sum_{\vy \in [0,R]^d } (1-f_{G}(\vy))
			+ \frac d{n^{10}}
			\leq {2(R+1)^d\over n}. \label{eq:EB_multi0}
		\end{align}
		where (a) followed as $\{\vY_i \in [0,R]^d\ \forall i=1,\dots,n\}$ implies
		$p^{\sfem}_n(\vy)=0$ for any $y\notin [0,R]^d$; and (b) follows from 
		$\Expect[p^\sfem_n(\vy)]=f_G(\vy)$ and, thus, $\EE[(p^\sfem_n(\vy) - f_G(\vy))^2] =
		\Var(p^\sfem_n(\vy))=\frac
		1{n^2}\sum_{i=1}^n\Var(\indc{\vY_i=\vy})={f_G(\vy)(1-f_G(\vy))\over n}$, and due to the choice of $R$ with \prettyref{lmm:support-multidim}. 		
		
		Using the union bound and the fact $H^2\leq 2$ we have 
		$$\EE\qth{H^2(f_{G},f_{{\hatGn}})\indc{\vY_i\notin [0,R]^d \text{ for some } i\in \{1,\dots,n\}}}
		\leq \frac {2d}{n^9}.$$
		Combining this with \prettyref{eq:EB_multi0} yields
		\begin{align*}
			\EE\qth{H^2(f_G,f_{{\hatGn}})}
			&\leq \EE\qth{H^2(f_{G},f_{{\hatGn}})\indc{\vY_i \in [0,R]^d\ \forall i=1,\dots,n}}
			+\EE\qth{H^2(f_{G},f_{{\hatGn}})\indc{\vY_i\notin [0,R]^d \text{ for some } i\in \{1,\dots,n\}}}
			\nonumber\\
			&\leq {4(R+1)^d\over n}
			\,,
		\end{align*}
		which completes the proof.
	\end{proof}
	
	\subsection{Regret bounds in multiple dimensions}
	
	\begin{proof}[Proof of \prettyref{thm:main_multidim}~]
		
	We first note that it suffices to only prove the case where the data generating $G$ satisfies $G\in \calP([0,h]^d)$. To prove the case where the marginals of $G$ belong to $\subexpo(s)$, it suffices to choose $h=\tilde c \max\{1,s\}\log n$ for a large enough constant $\tilde c>0$, as the following argument shows. Using the property of the Poisson mixture and the result on the support of the unconstrained NPMLE $\hat G$ for the one dimensional case in \prettyref{lmm:EB_support-g-hat}, as $f_{\vtheta}(\vy)=\prod_{j=1}^d f_{\theta_j}(y_j)$, we get
	\begin{enumerate}[label=\arabic*.~]
		\item $\hat G$ is supported on $[0,\max_{j=1}^d \max_{i=1}^nY_{ij}+1]^d$, which itself is a subset of $[0,\tilde c \max\{1,s\}\log n]^d$ with probability at least $1-\frac d{n^9}$ for a large enough constant $\tilde c>0$.
		\item As a result of the above and \prettyref{lmm:support-multidim}, we get that with probability at least $1-\frac d{n^9}$, each coordinate of $\hat \vtheta_{\hat G}$ lies in the interval $[0,\tilde c \max\{1,s\}\log n]$.
	\end{enumerate}
	Hence, using arguments similar to the one dimensional case in \eqref{eq:EB_m1} we can argue the following. 
	
	\begin{lemma}\label{lmm:regret_reduction_subexpo}
	For any $G$ with marginals in $\subexpo(s)$ and an estimate $\hat G$ supported on $[0,\hat h]^d$,
	$$
	\Regret(\hat G,G) \leq \Regret(\hat G,G_h) +O \pth{\hat h^2+s^2\over n^2},\quad h=\tilde c \max\{1,s\}\log n,
	$$
	where $G_h$ denote its restriction of $G$ on the hypercube $[0,h]^d$, i.e., $G_h[\vtheta\in \cdot ] = G[\vtheta\in \cdot|\vtheta\in [0,h]^d]$. 
	\end{lemma}
	We will prove the above result at the end of this section. In view of this, it suffices to bound $\Regret(\hat G,G_h)$ to get the desired regret upper bound.	To bound $\Regret(\hat G,G_h)=\EE_{G_h}\qth{\|\hat\vtheta_{\hat{G}}(\vY)-\hat\vtheta_{{G}_h}(\vY)\|^2}$ we use the following decomposition that is similar to the decomposition in the proof of the one dimensional case. Note that $\hat G$ is supported on $[0,\hat h]^d$, where $\hat h = h$ when $\hatGn$ is chosen to be supported over $[0,h]^d$ or as $\max_{j=1}^d \max_{i=1}^nY_{ij}+1$ when $\hatGn$ is obtained by performing an unconstrained optimization. Note that in the later case, as we argued above, $\hat h$ is bounded from above by $\tilde c s\log n$ with a probability $1-\frac d{n^9}$. Hence, we have
	\begin{align*}
		&\EE_{G_h}\qth{\|\hat\vtheta_{\hat{G}}(\vY)-\hat\vtheta_{{G}_h}(\vY)\|^2\indc{\vY\in [0,R]^d}}
		\nonumber\\
		&=\sum_{\vy\in [0,R]^d} \sum_{j=1}^d(y_j+1)^2f_{G_h}(\vy)\pth{{f_{\hat{G}}(\vy+\ve_j)\over f_{\hat{G}}(\vy)}-{f_{G_h}(\vy+\ve_j)\over f_{G_h}(\vy)}}^2
		\nonumber\\
		&\stepa{\leq} \sum_{\vy\in [0,R]^d} \sum_{j=1}^d (y_j+1)^2f_{G_h}(\vy)
		\Biggl\{3\pth{{f_{\hat{G}}(\vy+\ve_j)\over f_{\hat{G}}(\vy) }-{2f_{\hat{G}}(\vy+\ve_j)\over f_{G_h}(\vy)+f_{\hat{G}}(\vy)}}^2
		\nonumber\\
		&\quad +3\pth{{f_{G_h}(\vy+\ve_j)\over f_{G_h}(\vy) }-{2f_{G_h}(\vy+\ve_j)\over f_{G_h}(\vy)+f_{\hat{G}}(\vy)}}^2 +3\pth{2f_{G_h}(\vy+\ve_j)-2f_{\hat{G}}(\vy+\ve_j)\over f_{G_h}(\vy)+f_{\hat{G}}(\vy)}^2\Biggr\}
		\nonumber\\
		&\leq 3\sum_{\vy\in [0,R]^d} \sum_{j=1}^d\left\{
		\pth{(y_j+1)f_{\hat{G}}(\vy+\ve_j)\over f_{\hat{G}}(\vy) }^2{(f_{G_h}(\vy) -f_{\hat{G}}(\vy) )^2\over f_{G_h}(\vy)+f_{\hat{G}}(\vy)}\right.
		\nonumber\\
		&\quad \left. +\pth{(y_j+1)f_{G_h}(\vy+\ve_j)\over f_{G_h}(\vy)}^2{(f_{G_h}(\vy) -f_{\hat{G}}(\vy) )^2\over f_{G_h}(\vy)+f_{\hat{G}}(\vy)}
		+4(y_j+1)^2{(f_{G_h}(\vy+\ve_j)-f_{\hat{G}}(\vy+\ve_j))^2
			\over f_{G_h}(\vy)+f_{\hat{G}}(\vy)}\right\}
		\nonumber\\
		&\leq 3(h^2+\hat h^2)\sum_{\vy\in [0,R]^d} \sum_{j=1}^d{(f_{G_h}(\vy)-f_{\hat{G}}(\vy))^2\over f_{G_h}(\vy)+f_{\hat{G}}(\vy)}
		+12\sum_{\vy\in [0,R]^d} \sum_{j=1}^d(y_j+1)^2{(f_{G_h}(\vy+\ve_j)-f_{\hat{G}}(\vy+\ve_j))^2
			\over f_{G_h}(\vy)+f_{\hat{G}}(\vy)}
	\end{align*}
	
	where (a) followed from $(x+y+z)^2\leq 3(x^2+y^2+z^2)$ for any $x,y,z\in\reals$. Using
	$(\sqrt {f_{G_h}(\vx)}+\sqrt {f_{\hat{G}}(\vx)})^2\leq 2(f_{G_h}(\vx)+f_{\hat{G}}(\vx))$ for $\vx=\vy,\vy+\ve_j$ we continue the last display to get
	\begin{align*}
		&~\EE_{G_h}\qth{\|\hat\vtheta_{\hat{G}}(\vY)-\hat\vtheta_{{G}_h}(\vY)\|^2\indc{\vY\in [0,R]^d}}
		\nonumber\\
		&\leq 6(h^2+\hat h^2)\sum_{\vy\in [0,R]^d}\sum_{j=1}^d(\sqrt{f_{G_h}(\vy)}-\sqrt{f_{\hat{G}}(\vy)})^2
		\nonumber\\
		&\quad +24R\max_{\vy\in [0,R]^d}\sum_{j=1}^d{(y_j+1)f_{G_h}(\vy+\ve_j)+(y_j+1)f_{\hat{G}}(\vy+\ve_j)\over f_{G_h}(\vy)+f_{\hat{G}}(\vy)}
		\nonumber\\
		&\quad \cdot \sum_{\vy\in [0,R]^d}\sum_{j=1}^d(\sqrt{f_{G_h}(\vy+\ve_j)}-\sqrt{f_{\hat{G}}(\vy+\ve_j)})^2
		\nonumber\\
		&\leq d\pth{6(h^2+\hat h^2)+24(h+\hat h)R} {H^2(f_{\hat{G}},f_{G_h})},
	\end{align*}
	where the last inequality followed as
	\begin{align*}
	&\sum_{j=1}^d{(y_j+1)f_{G_h}(\vy+\ve_j)+(y_j+1)f_{\hat{G}}(\vy+\ve_j)\over f_{G_h}(\vy)+f_{\hat{G}}(\vy)}
	\sum_{\vy\in [0,R]^d}\sum_{j=1}^d(\sqrt{f_{G_h}(\vy+\ve_j)}-\sqrt{f_{\hat{G}}(\vy+\ve_j)})^2
	\nonumber\\
	&\leq \sum_{j=1}^d{(y_j+1)f_{G_h}(\vy+\ve_j)+(y_j+1)f_{\hat{G}}(\vy+\ve_j)\over f_{G_h}(\vy)+f_{\hat{G}}(\vy)}
	\sum_{\vy\in [0,R+1]^d}(\sqrt{f_{G_h}(\vy)}-\sqrt{f_{\hat{G}}(\vy)})^2
	\leq H^2(f_{\hat{G}},f_{G_h}).
	\end{align*}
	 Again, using $\hat\vtheta_{{G_h}}(\vy)\leq h,\hat\vtheta_{\hat {G}}(\vy)\leq \hat h$ we bound $\EE_{G_h}\qth{(\hat\vtheta_{{\hat G}}(\vY)-\hat\vtheta_{G_h}(\vY))^2\indc{\vY\notin [0,R]^d}}$ by $(h+\hat h)^2\epsilon_R(G_h)$, where $\epsilon_R(G_h)=\PP_{G_h}\qth{\vY\notin [0,R]^d}$. Combining this with the last display we get
	\begin{align*}
		\EE_{G_h}\qth{\|\hat\vtheta_{{\hat G}}(\vY)-\hat\vtheta_{G_h}(\vY)\|^2}\leq d\sth{6(h^2+\hat h^2)
			+24(h+\hat h)R} {H^2(f_{\hat{G}},f_{G_h})}+(h+\hat h)^2\epsilon_R(G_h).
	\end{align*}
	Finally we take expectation on both sides with respect to the training sample, $\hat G$, and $\hat h$. 
	\begin{enumerate}[label = (\roman*)~]
		\item  The proof of the result for $G\in\calP([0,h]^d)$, with both unconstrained and constrained estimator (which has knowledge of $h$) is very similar to the proof of \prettyref{thm:EB_main} and is omitted here.
		\item In the case of subexponential marginals of $G$, we pick $h=\tilde c \max\{1,s\} \log n$, for a large enough constant $\tilde c$. Hence, the support parameter $\hat h$ of $\hatGn$ as well as the support of each coordinate of $\vtheta$ is bounded from above by $h$ with a high probability. Using the high probability bound on $\hat h$ and the bound on $\EE\qth{{H^2(f_{\hat{G}},f_{G_h})}}$ as in \prettyref{thm:density_multidim}, the analysis for the bounded prior setup applies even though $\hatGn$ is obtained via an unconstrained optimization. Finally, using the bound on $\epsilon_R(G_h)$ as in \prettyref{lmm:support-multidim} we get the result.
	\end{enumerate}
	This finishes the proof.
\end{proof}

\begin{proof}[Proof of \prettyref{lmm:regret_reduction_subexpo}]
	We first note the multi-dimensional version of \eqref{eq:EB_m7} as follows:
	\begin{align*}
		\EE_{G}\qth{\|\hat\vtheta_{\hat G}(\vY)-\vtheta\|^2}
		&\leq 
		\EE_{G}\qth{(\hat\vtheta_{{{\hat G}}}(\vY)-\vtheta\|^2\indc{\vtheta\in [0,h]^d}}
		+\EE_{G}\qth{\|\hat\vtheta_{{{\hat G}}}(\vY)-\vtheta\|^2\indc{\vtheta\notin [0,h]^d}}
		\nonumber\\
		&\stepa{\leq}
		\EE_{G}\qth{\left.\|\hat\vtheta_{{{\hat G}}}(\vY)-\vtheta\|^2\right|{\vtheta\in [0,h]^d}}
		+\sqrt{\EE_{G}\qth{\|\hat\vtheta_{{{\hat G}}}(\vY)-\vtheta\|^4}\EE_G\qth{\indc{\vtheta\notin [0,h]^d}}}
		\nonumber\\
		&\stepb{\leq} \EE_{G_{h}}\qth{\|\hat\vtheta_{{{\hat G}}}(\vY)-\vtheta\|^2}+O\pth{\sqrt{d^4(\hat h^4+s^4){d\over n^9}}},
	\end{align*}
	where step (a) followed by Cauchy-Schwarz inequality, step (b) followed as 
	\begin{enumerate}[label=(\roman*)~]
		\item $(x+y)^4\leq 8(x^4+y^4)$ for any $x,y\in \reals$
		\item each coordinate of $\hat\vtheta_{{{\hat G}}}(\vY)$ is bounded by $\hat h$
		\item $\EE[\theta_j^4]\leq O(\max\{1,s^4\}),j=1,\dots,d$ by \prettyref{lmm:subexpo_properties}
		\item $\PP[\vtheta\notin [0,h]^d]\leq {d\over n^9}$ for a large enough $\tilde c$.
	\end{enumerate}
	Then, similar to \eqref{eq:EB_m1} in the one-dimensional case, the following equation applies (note that here $\hatGn$ is supported on $[0,\tilde c \max\{1,s\}\log n]$ with a high probability): 
	\begin{align*}\label{eq:subexp_truncate}
		\Regret(\hatGn; G)
		\le \Regret(\hatGn; G_h)
		+\mmse(G_h) - \mmse(G) + O\pth{\hat h^2+s^2\over n^2}.
	\end{align*}
	Note that we have $\mmse(G_h)\leq {\mmse(G)\over \PP\qth{\vtheta\in [0,h]^d}}$ from the following
	\begin{align*}
		\text{mmse}(G)
		&=\EE[\|\hat\vtheta_G - \vtheta\|^2]
		\geq \PP\qth{\vtheta\in [0,h]^d} \ \EE_{\vtheta\sim G}[\|\hat\vtheta_G - \vtheta\|^2 | \vtheta\in [0,h]^d]
		\ge \PP\qth{\vtheta\in [0,h]^d} \mmse(G_h).
	\end{align*}
	As $\text{mmse}(G)\le d$ (given that the naive estimation of $\vY$ achieves an expected loss of $d$) and $\PP\qth{\vtheta\notin [0,h]^d}\leq \frac 1{n^8}$, we get the desired result.
\end{proof}
	
\end{document}